\documentclass[12pt]{amsart}
\newtheorem{thm}{Theorem}[section]

\newtheorem{lem}[thm]{Lemma}
\newtheorem{prop}[thm]{Proposition}
\theoremstyle{definition}

\theoremstyle{remark}

\numberwithin{equation}{section}
\begin{document}

\begin{center}
{\bf{HYPERSURFACES IN SPACE FORMS\\
SATISFYING SOME GENERALIZED EINSTEIN METRIC CONDITION}}
\end{center}

\vspace{2mm}

\begin{center}
Ryszard Deszcz, Ma\l gorzata G\l ogowska and Georges Zafindratafa
\end{center}

\vspace{2mm}

\begin{center}
{\sl{Dedicated to Professor Leopold Verstraelen on his seventieth birthday}}
\end{center}

\vspace{2mm}

\noindent
{\bf{Abstract.}} 
The difference tensor $C \cdot R - R\cdot C$ 
of Einstein manifolds, some quasi-Einstein manifolds and Roter type manifolds, 
of dimension $n \geq 4$,
satisfy the following curvature condition: 
$(\ast )$ $C \cdot R - R\cdot C  = Q(S,C) - (\kappa/(n-1))\, Q(g,C)$.
We investigate hypersurfaces $M$ in space forms $N$ satisfying $(\ast )$.
The main result states that if the tensor $C \cdot R - R\cdot C$
of a non-quasi-Einstein hypersurface $M$ in $N$ is a linear combination of the tensors 
$Q(g,C)$ and $Q(S,C)$ then $(\ast )$ holds on $M$. 
In the case when $M$ is a quasi-Einstein hypersurface in $N$ and some additional assumptions are satisfied 
then $(\ast )$ also holds on $M$.

\section{Pseudosymmetry type curvature conditions}

Let $(M,g)$, $n = \dim M \geq 3$, be a semi-Riemannian manifold.
We denote by $\nabla$, $R$, $S$, $\kappa$ and $C$ 
the Levi-Civita connection, the Riemann-Christoffel curvature tensor, the Ricci tensor,
the scalar curvature and the Weyl conformal curvature tensor of $(M,g)$, respectively.
Let ${\mathcal U}_{S}$ be the set of all points of a manifold $(M,g)$, $n \geq 3$, at which 
$S$ is not proportional to $g$, i.e.,
${\mathcal U}_{S} \, = \,  \{x \in M\, | \, 
S - ( \kappa / n)\, g \neq 0\ \mbox {at}\ x \}$ and let
${\mathcal U}_{C}$ be the set of all points of a manifold $(M,g)$, $n \geq 4$, 
at which $C \neq 0$. We refer to sections 2 and 3 of the paper for precise definitions of the symbols used.

A semi-Riemannian manifold $(M,g)$, $n \geq 3$,
is said to be an {\sl Einstein manifold} \cite{Besse-1987} if at every point of $M$ 
its Ricci tensor $S$ is proportional to the metric tensor $g$, 
i.e., on $M$ we have
\begin{eqnarray}
S &=& \frac{\kappa }{n}\, g . 
\label{einstein000}
\end{eqnarray}
According to 
{\cite[p. 432] {Besse-1987}},
(\ref{einstein000}) is called the {\sl Einstein metric condition}.
Einstein manifolds form a natural subclass
of several classes of semi-Riemannian manifolds which are determined by
curvature conditions imposed on their Ricci tensor 
{\cite[Table, pp. 432-433] {Besse-1987}}.
These conditions are named {\sl generalized Einstein curvature conditions}
{\cite[Chapter XVI] {Besse-1987}}.
Any $3$-dimensional, as well as any conformally flat manifold 
$(M,g)$ of dimension $\geq 4$ satisfies 
\begin{eqnarray}
C \cdot R - R \cdot C  &=&   Q(S,C) - \frac{\kappa}{ n-1 }\, Q(g,C) ,
\label{Roterformula}
\end{eqnarray}
i.e. the condition $(\ast )$ mentioned in Abstract of this paper.  
It is easy to check that (\ref{Roterformula}) holds on any Einstein manifold (see Theorem 2.3).
Thus (\ref{Roterformula}) is a generalized Einstein metric condition. 
We refer to 
\cite{{P119}, {Ch-DDGP}, {DGHSaw}, {2003_DGHV}, {DGHZ01}, {2016_DGHZhyper}, 
{DGP-TV01}, {DH-2003}, {DHS-2001}, {DHS105}, {DP-TVZ}}
for results on semi-Riemannian manifolds, in particular, hypersurfaces in space forms, 
satisfying generalized Einstein metric conditions.

It is well-known that if a semi-Riemannian manifold $(M,g)$, $n \geq 3$, 
is locally symmetric then $\nabla R = 0$ on $M$
(see, e.g., {\cite[Chapter 1.5] {Lumiste}}). 
This implies the following integrability condition
${\mathcal{R}}(X,Y ) \cdot R = 0$, or briefly, 
\begin{eqnarray}
R \cdot R &=& 0 .
\label{semisymmetry}
\end{eqnarray}
Semi-Riemannian manifold satisfying (\ref{semisymmetry})
is called {\sl semisymmetric} (see, e.g., 
{\cite[Chapter 8.5.3] {TEC_PJR_2015}}, {\cite[Chapter 20.7] {Chen-2011}}, 
{\cite[Chapter 1.6] {Lumiste}}, \cite{{Szabo}, {LV3-Foreword}}).
Semisymmetric manifolds form a subclass of the class of pseudosymmetric manifolds.
A semi-Riemannian manifold $(M,g)$, $n \geq 3$, is said to be {\sl pseudosymmetric} 
if the tensors $R \cdot R$ and $Q(g,R)$ are linearly dependent at every point of $M$
(see, e.g., {\cite[Chapter 8.5.3] {TEC_PJR_2015}}, 
{\cite[Chapter 20.7] {Chen-2011}},
{\cite[Chapter 6] {DHV2008}},
{\cite[Chapter 12.4] {Lumiste}}, {\cite[Chapter 7.3.1] {MSV_2015}}, 
\cite{{D-1992}, {DGHHY}, {DGHSaw}, {DVV1991}, {HV_2007}, {HaVerSigma}, {KowSek_1997}, {SDHJK}, {LV1}, {LV2}, {LV3-Foreword}}). 
This is equivalent to
\begin{eqnarray}
R \cdot R &=& L_{R}\, Q(g,R) 
\label{pseudo}
\end{eqnarray}
on ${\mathcal{U}_{R} = \{x \in M\, | \, R - (\kappa / (n-1)n} )\, G \neq 0\ \mbox {at}\  x \}$,
where $L_{R}$ is some function on this set. 
According to \cite{KowSek_1997}, if the function $L_{R}$ is constant 
then $(M,g)$ is called a {\sl pseudosymmetric manifold of constant type}. 
Examples of non-semisymmetric pseudosymmetric manifolds are presented among others 
in \cite{{1989_DDV}, {DeKow}, {DVV1991}}. 
We also note that (\ref{pseudo}) implies
\begin{eqnarray}
(a)\  R \cdot S \ =\ L_{R}\, Q(g,S) 
\ \  &\mbox{and}&\ \ 
(b)\  R \cdot C \ =\ L_{R}\, Q(g,C) .
\label{Weyl-pseudo-bis}
\end{eqnarray}
The conditions 
(\ref{pseudo}), (\ref{Weyl-pseudo-bis})(a) and (\ref{Weyl-pseudo-bis})(b)
are equivalent on the set ${\mathcal{U}}_{S} \cap {\mathcal{U}}_{C}$ 
of any $4$-dimensional warped product manifold 
$\overline{M} \times _{F} \widetilde{N}$,
(see, e.g., \cite{DGJZ-2016} and references therein).
We note that
on any semi-Riemannian manifold $(M,g)$, $n \geq 3$, we have
${\mathcal{U}}_{S} \cup {\mathcal{U}}_{C} = {\mathcal{U}}_{R}$ (see, e.g., \cite{DGHHY}). 

A semi-Riemannian manifold $(M,g)$, $n \geq 3$, is called {\sl Ricci-pseudosymmetric} 
if the tensors $R \cdot S$ and $Q(g,S)$ are linearly dependent at every point of $M$
(see, e.g., {\cite[Chapter 8.5.3] {TEC_PJR_2015}}, \cite{{D-1992}, {DGHSaw}, {LV3-Foreword}}).
This is equivalent on ${\mathcal{U}}_{S} \subset M$ to 
\begin{eqnarray}
R \cdot S &=& L_{S}\, Q(g,S) , 
\label{Riccipseudo07}
\end{eqnarray}
where $L_{S}$ is some function on this set. 
According to \cite{G6}, if the function $L_{S}$ is constant 
then $(M,g)$ is called a {\sl Ricci-pseudosymmetric manifold of constant type}. 
Every warped product manifold $\overline{M} \times _{F} \widetilde{N}$
with an $1$-dimensional $(\overline{M}, \overline{g})$ manifold and
an $(n-1)$-dimensional Einstein semi-Riemannian manifold $(\widetilde{N}, \widetilde{g})$, $n \geq 3$, 
and a warping function $F$, 
is a Ricci-pseudosymmetric manifold
(see, e.g., {\cite[Chapter 8.5.3] {TEC_PJR_2015}}, {\cite[Section 1] {Ch-DDGP}}, {\cite[Example 4.1] {DGJZ-2016}}).
A special subclass of Ricci-pseudo\-sym\-metric manifolds form Ricci-semisymmetric manifolds.
A semi-Riemannian manifold $(M,g)$, $n \geq 3$
is said to be {\sl Ricci-semisymmetric} if 
\begin{eqnarray}
R \cdot S &=& 0
\label{Riccisemisymmetric}
\end{eqnarray}
on $M$ (see, e.g., {\cite[Chapter 8.5.3] {TEC_PJR_2015}}, \cite{DGHSaw}). 
Ricci-semisymmetric manifolds were investigated by several authors,
see, e.g., 
\cite{{AD-2002}, 
{TEC_PJR_2015},
{2002-DG-1}, {DG90},
{Lumiste}, 
{Mir-1992}, {Mir-Mach-2012}, {Ryan-1972}}
and references therein.
Ricci-semisymmetric manifolds (submanifolds) are
also named {\sl Ric-semisymmetric manifolds} (submanifolds), 
and in particular, {\sl Ric-semisymmetric hypersurfaces}
{\cite[Chapter 12.7] {Lumiste}}, \cite{{Mir-1992}, {Mir-Mach-2012}}
or {\sl Ryan hypersurfaces} {\cite[Chapter 8.5.3] {TEC_PJR_2015}}.

A semi-Riemannian manifold $(M,g)$, $n \geq 4$, is said to be {\sl Weyl-pseudo\-sym\-met\-ric} 
if the tensors $R \cdot C$ and $Q(g,C)$ are linearly dependent at every point of $M$
\cite{{D-1992}, {DGHHY}, {DGHSaw}}. 
This is equivalent on ${\mathcal{U}}_{C}\subset M$ to 
\begin{eqnarray}
R \cdot C &=& L_{1}\, Q(g,C) ,  
\label{Weyl-pseudo}
\end{eqnarray}
where $L_{1}$ is some function on this set. 
In particular, if the condition $R \cdot C = 0$ holds on ${\mathcal{U}}_{C} \subset M$
then $(M,g)$, $n \geq 4$, is called {\sl Weyl-semisymmetric}
{\cite[Chapter 20.7] {Chen-2011}}, \cite{{DGHSaw}, {DHV2008}, {LV1}}.

A semi-Riemannian manifold $(M,g)$, $n \geq 4$, is said to have {\sl pseudosymmetric Weyl tensor}
if the tensors $C \cdot C$ and $Q(g,C)$ are linearly dependent at every point of $M$ 
(see, e.g., {\cite[Chapter 20.7] {Chen-2011}}, \cite{{DGHHY}, {DGHSaw}, {DGJZ-2016}}).
This is equivalent on ${\mathcal U}_{C}\subset M$ to 
\begin{eqnarray}
C \cdot C &=& L_{C}\, Q(g,C) ,  
\label{4.3.012}
\end{eqnarray}
where $L_{C}$ is some function on this set. 
Every warped product manifold
$\overline{M} \times _{F} \widetilde{N}$,
with $\dim \overline{M}  = \dim \widetilde{N} = 2$, 
satisfies (\ref{4.3.012})
(see, e.g., \cite{{DGHHY}, {DGHSaw}, {DGJZ-2016}} and references therein).
Thus in particular,
the Schwarzschild spacetime, the Kottler spacetime
and the Reissner-Nordstr\"{o}m spacetime satisfy (\ref{4.3.012}).
Recently, manifolds satisfying (\ref{4.3.012})
were investigated among others in \cite{{DGHHY}, {DGJZ-2016}, {DeHoJJKunSh}}.

Warped product manifolds $\overline{M} \times _{F} \widetilde{N}$, of dimension $\geq 4$,
satisfying on 
${\mathcal U}_{C} \subset \overline{M} \times _{F} \widetilde{N}$ the condition
\begin{eqnarray}
R \cdot R - Q(S,R) &=& L\, Q(g,C) ,  
\label{genpseudo01}
\end{eqnarray}
where $L$ is some function on this set,
were studied among others in \cite{{49}, {DGJZ-2016}}.
For instance, 
in \cite{49} necessary and sufficient conditions for  
$\overline{M} \times _{F} \widetilde{N}$ to be a manifold satisfying (\ref{genpseudo01}) are given.
Moreover, in that paper it was proved that 
any $4$-dimensional warped product manifold $\overline{M} \times _{F} \widetilde{N}$, 
with an $1$-dimensional base $(\overline{M},\overline{g})$, 
satisfies (\ref{genpseudo01}) {\cite[Theorem 4.1] {49}}.
The warped product manifold $\overline{M} \times _{F} \widetilde{N}$, 
with $2$-dimensional base $(\overline{M},\overline{g})$ 
and $(n-2)$-dimensional space of constant curvature $(\widetilde{N},\widetilde{g})$, $n \geq 4$,
is a manifold satisfying 
(\ref{4.3.012}) and (\ref{genpseudo01}) {\cite[Theorem 7.1 (i)] {DGJZ-2016}}.

We refer to
\cite{{Ch-DDGP}, {D-1992}, {DGHHY}, {DGHSaw}, {DGHZ01}, {DGJZ-2016}, {DHV2008}, {DeHoJJKunSh}, {DP-TVZ}, {SDHJK}, {LV1}} 
for results on semi-Riemannian manifolds satisfying (\ref{pseudo})-(\ref{genpseudo01}), 
as well as other conditions of this kind, named 
{\sl{pseudosymmetry type curvature conditions}}
or 
{\sl{pseudosymmetry type conditions}}. 
It seems that (\ref{pseudo}) 
is the most important condition of that family of curvature conditions (see, e.g., \cite{DGJZ-2016}).
We also can state that
the Schwarzschild spacetime, the Kottler spacetime, the Reissner-Nordstr\"{o}m spacetime, 
as well as the Friedmann-Lema{\^{\i}}tre-Robertson-Walker spacetimes are the "oldest" examples 
of pseudosymmetric warped product manifolds (see, e.g., \cite{{DGJZ-2016}, {DHV2008}, {DVV1991}, {SDHJK}}).

A semi-Riemannian manifold $(M,g)$, $n \geq 3$, is said to be 
a {\sl quasi-Einstein manifold} if 
\begin{eqnarray}
\mathrm{rank}\, (S - \alpha\, g) &=& 1
\label{quasi02}
\end{eqnarray}
on ${\mathcal U}_{S} \subset M$, where $\alpha $ is some function on this set.
In Section 2 (see Remark 2.5 (i)-(iii)) we present some facts related to that class of manifolds. 
We mention that some quasi-Einstein warped product manifolds satisfy the following condition 
(see Remark 2.5 (iii))
\begin{eqnarray}
(n-2) ( R \cdot C - C \cdot R ) &=& Q(S,C) - L_{S}\, Q(g,C) .
\label{quasi022}
\end{eqnarray}
Quasi-Einstein manifolds arose during the study of exact solutions
of the Einstein field equations and the investigation on quasi-umbilical hypersurfaces 
of conformally flat spaces (see, e.g., \cite{{DGHSaw}, {DGJZ-2016}} and references therein). 
Recently quasi-Einstein manifolds satisfying some pseudosymmetry type conditions  
were investigated among others in
\cite{{P119}, {Ch-DDGP}, {DGHHY}, {DGHZ01}, {DeHoJJKunSh}}. 
Quasi-Einstein hypersurfaces in semi-Riemannian spaces of constant curvature
were studied among others in
\cite{{DGHS}, {R102}, {DHS105}, {G6}}, see also {\cite[Chapter 6.2] {TEC_PJR_2015}},
{\cite[Chapter 19.5] {Chen-2011}},
{\cite[Chapter 4.6] {Chen-2017}},
\cite{{DGHSaw}, {LV3-Foreword}} and references therein. 
We mention that there are different extensions of the class of quasi-Einstein manifolds. 
For instance we have the class of almost quasi-Einstein manifolds \cite{Chen-2017-KJM}
and the class of $2$-quasi-Einstein manifolds (see, e.g. \cite{{2016_DGHZhyper}, {DGJZ-2016}}).

Investigations on semi-Riemannian manifolds $(M,g)$, $n \geq 4$,
satisfying 
(\ref{pseudo}) and (\ref{4.3.012})
or
(\ref{pseudo}) and (\ref{genpseudo01})
on ${\mathcal{U}}_{S} \cap {\mathcal{U}}_{C} \subset M$
lead to the following condition ({\cite[Theorem 3.2 (ii)] {DY1994}}, {\cite[Lemma 4.1] {P43}}, see also
{\cite[Section 1] {DGJZ-2016}})
\begin{eqnarray}
R &=& \frac{\phi}{2}\, S\wedge S + \mu\, g\wedge S + \frac{\eta}{2}\, g \wedge g ,
\label{eq:h7a}
\end{eqnarray}
where 
$\phi$, $\mu $ and $\eta $ are some functions on ${\mathcal U}_{S} \cap {\mathcal U}_{C}$.
A semi-Riemannian manifold $(M,g)$, $n \geq 4$, satisfying (\ref{eq:h7a}) on 
${\mathcal U}_{S} \cap {\mathcal U}_{C} \subset M$ 
is called a {\sl Roter type manifold}, or a {\sl Roter type space}, or a {\sl Roter space} 
\cite{{P106}, {DGP-TV01}, {DGP-TV02}}. 
Roter type manifolds and in particular Roter type hypersurfaces
(i.e. hypersurfaces satisfying (\ref{eq:h7a})), 
in semi-Riemannian spaces of constant curvature were studied in:
\cite{{P106}, {DGHHY}, {DGHZ01}, {DGP-TV01}, {R102}, {DeKow}, {DP-TVZ}, {DePlaScher}, {DeScher}, {G5}, {Kow01}, {Kow02}}. 
Roter type manifolds satisfy (\ref{Roterformula}), 
as well as some other pseudosymmetry type conditions (see Theorem 2.4 and Remark 2.5 (iv)-(vii)). 
We note that if (\ref{eq:h7a}) is satisfied 
at a point of ${\mathcal U}_{S} \cap {\mathcal U}_{C}$ then at this point 
\begin{eqnarray} 
\mathrm{rank} ( S - \alpha \, g ) > 1 ,\ \ \mbox{for any}\ \ \alpha \in  {\mathbb{R}} .
\label{dhs_quasi_Einstein}
\end{eqnarray}

Let $M$, $n \geq 3$, be a connected hypersurface 
isometrically immersed in a semi-Riemannian space of constant curvature 
$N_{s}^{n+1}(c)$.
Let 
$g$ be the metric tensor induced on $M$ from the metric of the ambient space 
and let $R$, $S$, $\kappa$ and $C$
be the Riemann-Christoffel curvature tensor, Ricci tensor, 
the scalar curvature and the Weyl conformal curvature tensor of $g$, respectively. 
Further, let $H$ and ${\mathcal A}$ be the second fundamental tensor and the shape operator
of $M$, respectively. 
We have $H(X,Y) = g( {\mathcal A}X,Y)$, for any vector fields $X,Y$ tangent to $M$.
The $(0,2)$-tensors $H^{2}$ and $H^{3}$ are defined by
$H^{2}(X,Y) = H( {\mathcal A}X,Y)$ and $H^{3}(X,Y) = H^{2}( {\mathcal A}X,Y)$, respectively.
Hypersurfaces in $N_{s}^{n+1}(c)$ satisfying pseudosymmetry type conditions were investigated in 
several papers, see, e.g., {\cite[Section 1] {Saw-2005}} and references therein.
We also refer to \cite{{2016_DGHZhyper}, {DGP-TV02}, {Saw-2015}}
for recent results related to this subject. For instance,
if $M$ is a hypersurface in $N_{s}^{n+1}(c)$, $n \geq 3$, satisfying 
\begin{eqnarray}
H^{2} &=& \alpha \, H + \beta \, g ,
\label{two}
\end{eqnarray}
for some functions $\alpha$ and $\beta$ on $M$, 
then $M$ is a pseudosymmetric manifold (see Theorem 3.1). 
If some additional assumptions are satisfied then $M$ is a Roter type manifold (see Theorem 3.3).  
In particular, every non-Einstein and non-conformally flat Clifford torus of dimension $\geq 5$
is a Roter type manifold (see Example 3.6). 

Let $M$ be a hypersurface in $N_{s}^{n+1}(c)$, $n \geq 4$. 
We denote by ${\mathcal U}_{H} \subset M$ the set 
of all points at which the tensor $H^{2}$ 
is not a linear combination of the tensors $H$ and $g$ of $M$.
It is known that
${\mathcal U}_{H} \subset {\mathcal U}_{S} \cap {\mathcal U}_{C} \subset M$
({\cite[Proposition 2.1] {ADEMO-2002}}, {\cite[Section 1] {P106}}). 
For instance, if $M$ is the Cartan hypersurface then ${\mathcal U}_{H} = M$
({\cite[Chapter 3.8.3] {TEC_PJR_2015}}, {\cite[Chapter 20.3] {Chen-2011}}).
In addition, if $M$ is the Cartan hypersurface of dimension $6$, $12$, or $24$ then it is 
a non-pseudosymmetric Ricci-pseudosymmetric manifold of constant type satisfying
{\cite[Proposition 1 (i), Theorem 1] {DY-1994-cm}}
\begin{eqnarray}
R \cdot S &=& \frac{\widetilde{\kappa }}{ n (n+1) }\, Q(g,S) ,
\label{DZ004Ricciaaa}
\end{eqnarray}
provided that $n = 6, 12, 24$.
The $3$-dimensional Cartan hypersurface $M$ is a pseudosymmetric manifold of constant type. 
In fact, 
$R \cdot R = ( \widetilde{\kappa}/12 )\, Q(g,R)$ on $M$ 
{\cite[Section 5] {DVY}}. 

Evidently, (\ref{semisymmetry}) implies (\ref{Riccisemisymmetric}).
The converse is not true. 
The problem of the equivalence of (\ref{semisymmetry}) and (\ref{Riccisemisymmetric}) on hypersurfaces,
named {\sl Ryan's problem} {\cite[Chapter 12.7] {Lumiste}},
was investigated by several authors 
(see, e.g., {\cite[Chapter 8.5.3] {TEC_PJR_2015}},
\cite{{AD-2002}, {2000_DDKV}, {1997_DDSVY}, {1999_DDSVY}, {Mir-1999}} and references therein). 
This problem was stated as Problem P808 in \cite{Ryan-1972} (cf., {\cite[Chapter 12.7] {Lumiste}}).
We mention that (\ref{semisymmetry}) and (\ref{Riccisemisymmetric})
are equivalent on hypersurfaces in $5$-dimensional semi-Riemannian spaces of
constant curvature $N_{s}^{5}(c)$ \cite{1999_DDSVY}. 
For a presentation of results on the problem of the equivalence 
of semisymmetry, Ricci-semisymmetry and Weyl-semisymmetry, 
or, more generally, 
of pseudosymmetry, Ricci-pseudo\-sym\-met\-ry and Weyl-pseudosymmetry
on semi-Riemannian manifolds, and, in particular, 
on hypersurfaces in semi-Riemannian spaces,
we refer to {\cite[Section 4] {DGHSaw}}.  

The second fundamental tensor $H$ of hypersurfaces $M$ in $N_{s}^{n+1}(c)$, $n \geq 4$,
realizing some pseudosymmetry type conditions on ${\mathcal U}_{H} \subset M$
satisfy also the following equation
\begin{eqnarray}
H^{3} &=& tr(H)\, H^{2} + \psi \, H + \rho \, g ,
\label{DS4}
\end{eqnarray}
where $\psi $ and $\rho $ are some functions on ${\mathcal U}_{H}$.
We refer to 
\cite{{2003_DGHV}, {2016_DGHZhyper}, {2011-DGPSS}, {Saw-2004}, {Saw-2005}, {Saw-2006}, {Saw-2007}, {Saw-2015}}
for results on hypersurfaces satisfying (\ref{DS4}). 
Evidently, if $M$ is a hypersurface in an $4$-dimensional Riemannian space of constant curvature $N^{4}(c)$ 
then (\ref{DS4}) holds ${\mathcal U}_{H} \subset M$.
If $M$ is a hypersurface in $N_{s}^{5}(c)$ satisfying (\ref{DS4}) on ${\mathcal U}_{H} \subset M$
then we have ({\cite[Proposition 2.1] {2005-DSaw}}, Proposition 4.3 (i))
\begin{eqnarray}
H^{3} &=& tr(H)\, H^{2} + \psi \, H .
\label{DS4aa}
\end{eqnarray}

If (\ref{DS4aa}) is satisfied 
at every point of the set ${\mathcal U}_{H}$ of a hypersurface $M$
in  $N_{s}^{n+1}(c)$, $n \geq 3$, then (\ref{DZ004Ricciaaa}) holds on this set (see Proposition 4.3 (v)).
Further,
if $\mathrm{rank} H = 2$ at every point of the set ${\mathcal U}_{H}$ of a hypersurface $M$
in  $N_{s}^{n+1}(c)$, $n \geq 3$, then (\ref{DS4aa}) holds on this set (see Proposition 4.3 (ii), or   
{\cite[Lemma 2.1] {P106}}). In addition, on this set we have (see Theorem 3.2 (iii))
\begin{eqnarray}
R \cdot R &=& \frac{ \widetilde{\kappa} }{n(n+1)}\, Q(g,R) .
\label{pseudo-constant01}
\end{eqnarray}
If $\mathrm{rank} H = 2$ at every point of a hypersurface $M$ in $N_{s}^{n+1}(c)$, $n \geq 3$,  
then $M$ is called a {\sl hypersurface with type number} $2$ (see, e.g. \cite{Chen-Yildirim-2015}). Thus we see that if 
$M$ in $N_{s}^{n+1}(c)$, $n \geq 3$, 
is a hypersurface with type number $2$ then $M$ is a pseudosymmetric manifold of constant type.
We refer to \cite{Chen-Yildirim-2015} (see also {\cite[Section 5]{R102}}, \cite{DVZ-2018})
for examples of submanifolds, 
and in particular, hypersurfaces in spaces of constant curvature with type number $\leq 2$.
 
Evidently, (\ref{DS4}) is a particular case of the equation
\begin{eqnarray}
H^{3} &=& \phi\, H^{2} + \psi \, H + \rho \, g ,
\label{DS4_general}
\end{eqnarray}
where $\phi$, $\psi $ and $\rho $ are some functions on ${\mathcal U}_{H}$.
Hypersurfaces $M$ in $N_{s}^{n+1}(c)$, $n \geq 4$, 
satisfying (\ref{DS4_general}) on ${\mathcal U}_{H} \subset M$
were investigated among others in
\cite{{TEC_GRJ_1998}, {DGP-TV02}, {Saw-2004}, {Saw-2015}},
see also 
{\cite[Chapters 3.8.3 and 5.6]{TEC_PJR_2015}}
and references therein.
If the tensor $R \cdot C$, $C \cdot R$ or $R \cdot C - R \cdot C$ 
is on ${\mathcal U}_{H} \subset M$
a linear combination of the tensor $R \cdot R$ 
and of a finite sum of the Tachibana tensors of the form $Q(A,B)$, 
where $A$ is a symmetric $(0,2)$-tensor and $B$ a generalized curvature tensor,
then the tensor $H$ satisfies (\ref{DS4}) on this set
({\cite[Corollary 4.1] {2003_DGHV}}).  

In Proposition 4.1 we present results 
on hypersurfaces $M$ in $N_{s}^{n+1}(c)$, $n \geq 4$, satisfying (\ref{DS4})
obtained in {\cite[Proposition 5.1]{Saw-2005}}. 
Further, in that section we prove the following results.

Let $M$ be a hypersurface in $N_{s}^{n+1}(c)$, $n \geq 4$,
satisfying (\ref{DS4}) on ${\mathcal U}_{H} \subset M$. 
We have: 

(i) (see Proposition 4.2 (i)) The conditions  
\begin{eqnarray}
R \cdot S &=&  
Q(g, S^{2}) + \left( \varepsilon \psi - \frac{( 2 n - 3) \widetilde{\kappa }}{n(n+1)} \right) Q(g,S ) ,
\label{GGG01}\\
R \cdot S^{2} &=& Q(S, S^{2}) 
+ \rho _{1}\, Q(g, S^{2}) + \rho _{2}\, Q(g, S ) ,
\label{EEE01}\\
S^{3} &=& \left( - 2 \varepsilon \psi + \frac{ 3 (n - 1) \widetilde{\kappa }}{n(n+1)} \right) S^{2} + \rho_{2} \, S + \rho _{3}\, g ,
\label{EEE01new}
\end{eqnarray}
hold on ${\mathcal U}_{H}$, where $\rho _{1}$, $\rho _{2}$ and $\rho _{3}$ are defined by (\ref{EEE01bb}).

(ii) (see Theorem 4.6) 
If on ${\mathcal U}_{H}$
the tensor $Q(S,R)$ is equal to the Tachibana tensor $Q(g,T)$, 
where $T$ is a generalized curvature tensor, then any of the tensors:  
$R \cdot R$, $R \cdot C$, $C \cdot R$, $R \cdot C - C \cdot R$ and $C \cdot C$
is equal to some Tachibana tensor $Q(g,B)$, 
where $B$ is a linear combination of the tensors
$R$, $g \wedge g$, $g \wedge S$, $g \wedge S^{2}$ and $S \wedge S$.

(iii) (see Proposition 4.7) The following conditions are satisfied on ${\mathcal U}_{H}$
\begin{eqnarray}
C \cdot C 
&=&
\frac{n-3}{n-2}\, R \cdot C
+ \frac{1}{n-2} \left( \frac{\kappa }{n-1} + \varepsilon \psi
- \frac{( 2 n - 3) \widetilde{\kappa }}{n(n+1)} \right)
Q(g,C),
\label{DS16Anew01}\\
\nonumber\\
(n-2)\, C \cdot R  + R \cdot C 
&=& 
(n-2)\, Q(S,C) 
+ 
\left(
\frac{\kappa}{n-1} + \varepsilon \psi - \frac{ (n-1)^{2} \widetilde{\kappa}}{n (n+1)} 
\right)
Q(g,C)\nonumber\\
& & - \frac{1}{(n-2)}\, Q\left( g, \frac{n-2}{2}\, S \wedge S - \kappa \, g\wedge S + g \wedge S^{2} \right) .
\label{02identity01hyper17}
\end{eqnarray}

(iv) (see Theorem 4.8) If (\ref{quasi02}) and (\ref{DS4aa}) are satisfied on ${\mathcal U}_{H}$ 
then (\ref{Roterformula}) holds on ${\mathcal U}_{H}$  if and only if the following two conditions hold
on this set 
\begin{eqnarray}
(a)\ \ 
\frac{\kappa }{n-1} \ =\ \frac{\widetilde{ \kappa }}{n+1}
\ \ \ 
&\mbox{and}&
\ \ \ 
(b)\ \  Q\left( S - \frac{\kappa }{n}\, g, C \right)\ =\ 0 .
\label{qqee08}
\end{eqnarray}

In Section 5 we consider hypersurfaces $M$ in  $N_{s}^{n+1}(c)$, $n \geq 4$, 
satisfying on ${\mathcal U}_{H} \subset M$ 
\begin{eqnarray}
R \cdot C - C \cdot R &=& L_{1}\, Q(S,C) + L_{2}\, Q(g,C) ,
\label{cond01}
\end{eqnarray}
where $L_{1}$ and $L_{2}$ are some functions defined on this set.
Theorem 5.1 states that if on ${\mathcal U}_{H}$
the conditions (\ref{quasi02}) and (\ref{cond01}) are satisfied, 
for some functions $\alpha$, $L_{1}$ and $L_{2}$,
then on this set we have (\ref{DS4aa}), for some function $\psi$, and 
\begin{eqnarray}
(n-2)\, ( R \cdot C - C \cdot R ) &=& Q(S,C) 
- \frac{\widetilde{\kappa }}{ n (n+1) }\, Q(g,C) .
\label{quasi321}
\end{eqnarray}
Finally, theorems 5.2 and 5.3 state that if on ${\mathcal U}_{H}$
the conditions (\ref{dhs_quasi_Einstein}) and (\ref{cond01}) are satisfied,
for some functions $L_{1}$ and $L_{2}$, 
then on this set we have (\ref{Roterformula}) and $Q(S,R) = Q(g,T)$, i.e. (\ref{cond02uuu}), 
where $T$ is a linear combination of the tensors  
$R$, $g \wedge g$, $g \wedge S$, $g \wedge S^{2}$ and $S \wedge S$. 
Moreover, any of the tensors: $R \cdot R$, $R \cdot C$, $C \cdot R$, $R \cdot C - C \cdot R$ and $C \cdot C$
is equal to some Tachibana tensor $Q(g,B)$, 
where $B$ is a linear combination of the tensors
$R$, $g \wedge g$, $g \wedge S$, $g \wedge S^{2}$ and $S \wedge S$.

\section{Preliminary results}

Throughout this paper all manifolds are assumed
to be connected paracompact manifolds of class $C^{\infty }$.
Let $(M,g)$ be an $n$-dimensional
semi-Riemannian manifold
and let $\nabla$ be its
Levi-Civita connection
and $\Xi (M)$
the Lie algebra of vector fields on $M$.
We define on $M$ the endomorphisms
$X \wedge _{A} Y$ and
${\mathcal{R}}(X,Y)$
of $\Xi (M)$ by 
$(X \wedge _{A} Y)Z = A(Y,Z)X - A(X,Z)Y$ and
\begin{eqnarray*}
{\mathcal R}(X,Y)Z
&=&
\nabla _X \nabla _Y Z - \nabla _Y \nabla _X Z - \nabla _{[X,Y]}Z ,
\end{eqnarray*}
respectively, where $A$ is a symmetric $(0,2)$-tensor on $M$
and $X, Y, Z \in \Xi (M) $.
The Ricci tensor $S$,
the Ricci operator ${\mathcal{S}}$,
the tensors $S^{2}$ and $S^{3}$ and
the scalar curvature $\kappa $
of $(M,g)$ are defined by
$S(X,Y) = \mathrm{tr} \{ Z \rightarrow {\mathcal{R}}(Z,X)Y \}$,
$g({\mathcal{S}}X,Y) = S(X,Y)$,
$S^{2}(X,Y) = S({\mathcal{S}}X,Y)$,
$S^{3}(X,Y) = S^{2}({\mathcal{S}}X,Y)$
and
$\kappa \, =\, \mathrm{tr}\, {\mathcal{S}}$, respectively.
The endomorphism ${\mathcal{C}}(X,Y)$ of $(M,g)$, $n \geq 3$, is defined by
\begin{eqnarray*}
{\mathcal{C}}(X,Y)Z  &=& {\mathcal{R}}(X,Y)Z
- \frac{1}{n-2} \left( X \wedge _{g} {\mathcal{S}}Y + {\mathcal{S}}X \wedge _{g} Y
- \frac{\kappa}{n-1}X \wedge _{g} Y \right) Z .
\end{eqnarray*}
The $(0,4)$-tensor $G$,
the Riemann-Christoffel curvature tensor $R$ and
the Weyl conformal curvature tensor $C$ of $(M,g)$ are defined by
$G(X_1,X_2,X_3,X_4) = g (( (X_1 \wedge _{g} X_2)X_3,X_4)$,  
\begin{eqnarray*}
R(X_1,X_2,X_3,X_4) \ =\ g({\mathcal{R}}(X_1,X_2)X_3,X_4) ,\ \ \
C(X_1,X_2,X_3,X_4) \ =\ g({\mathcal{C}}(X_1,X_2)X_3,X_4) ,
\end{eqnarray*}
respectively, where $X_1,X_2,X_3,X_4 \in \Xi (M)$.
Let ${\mathcal{B}}$ be a tensor field sending any $X, Y \in \Xi (M)$
to a skew-symmetric endomorphism ${\mathcal{B}}(X,Y)$
and let $B$ be a $(0,4)$-tensor associated with ${\mathcal{B}}$ by
\begin{eqnarray}
B(X_1,X_2,X_3,X_4) &=&
g({\mathcal{B}}(X_1,X_2)X_3,X_4) .
\label{DS5}
\end{eqnarray}
The tensor $B$ is said to be a {\sl{generalized curvature tensor}}
if the following conditions are satisfied
\begin{eqnarray*}
& &
B(X_1,X_2,X_3,X_4) \ =\  B(X_3,X_4,X_1,X_2) ,\\
& &
B(X_1,X_2,X_3,X_4)
+ B(X_3,X_1,X_2,X_4)
+ B(X_2,X_3,X_1,X_4) \ =\  0 .
\end{eqnarray*}
For ${\mathcal{B}}$ as above, let $B$ be again defined by (\ref{DS5}).
We extend the endomorphism
${\mathcal{B}}(X,Y)$ to a derivation
${\mathcal{B}}(X,Y) \cdot \, $
of the algebra of tensor fields on $M$,
assuming that it commutes with contractions and
$\ {\mathcal{B}}(X,Y) \cdot \! f \, =\, 0$, for any smooth function $f$ on $M$.
For a $(0,k)$-tensor field $T$, $k \geq 1$,
we can define the $(0,k+2)$-tensor $B \cdot T$ by
\begin{eqnarray*}
& & (B \cdot T)(X_1,\ldots ,X_k,X,Y) \ =\
({\mathcal{B}}(X,Y) \cdot T)(X_1,\ldots ,X_k)\\
&=& - T({\mathcal{B}}(X,Y)X_1,X_2,\ldots ,X_k)
- \cdots - T(X_1,\ldots ,X_{k-1},{\mathcal{B}}(X,Y)X_k) .
\end{eqnarray*}
If $A$ is a symmetric $(0,2)$-tensor then we define
the $(0,k+2)$-tensor $Q(A,T)$ by
\begin{eqnarray*}
& & Q(A,T)(X_1, \ldots , X_k, X,Y) \ =\
(X \wedge _{A} Y \cdot T)(X_1,\ldots ,X_k)\\
&=&- T((X \wedge _A Y)X_1,X_2,\ldots ,X_k)
- \cdots - T(X_1,\ldots ,X_{k-1},(X \wedge _A Y)X_k) .
\end{eqnarray*}
The tensor $Q(A,T)$
is called the {\sl Tachibana tensor of the tensors} $A$ and $T$,
in short the Tachibana tensor 
(see, e.g., \cite{{DGHHY}, {DGHZ01}, {DGJZ-2016}, {2011-DGPSS}, {DHV2008}, {DeHoJJKunSh}}).
Thus, among other things, we have the $(0,6)$-tensors:
$R \cdot R$, $R \cdot C$, $C \cdot R$, $C \cdot C$, 
$Q(g,R)$, $Q(S,R)$, $Q(g,C)$ and $Q(S,C)$, as well as the $(0,4)$-tensors:
$R \cdot S$, $C \cdot S$, $Q(g,S)$, $Q(g,S^{2})$ and $Q(S,S^{2})$.  
For a symmetric $(0,2)$-tensor $E$ and a $(0,k)$-tensor $T$, $k \geq 2$,
we define their Kulkarni-Nomizu product $E \wedge T$ by 
(see, e.g., \cite{{DGHHY}, {2011-DGPSS}})
\begin{eqnarray*}
& &(E \wedge T )(X_{1}, \ldots , X_{4}, Y_{3}, \ldots , Y_{k})\\
&=&
E(X_{1},X_{4}) T(X_{2},X_{3}, Y_{3}, \ldots , Y_{k})
+ E(X_{2},X_{3}) T(X_{1},X_{4}, Y_{3}, \ldots , Y_{k} )\\
& &
- E(X_{1},X_{3}) T(X_{2},X_{4}, Y_{3}, \ldots , Y_{k})
- E(X_{2},X_{4}) T(X_{1},X_{3}, Y_{3}, \ldots , Y_{k}) .
\end{eqnarray*}

On any semi-Riemannian manifold $(M,g)$, $n \geq 3$, we have
(see, e.g., \cite{DGJZ-2016})
\begin{eqnarray}
C &=& R - \frac{1}{n-2}\, g\wedge S + \frac{\kappa }{(n-2)(n-1)}\, G ,  
\label{eqn2.1}\\
(a)\ \ Q(A,G)
&=&
- Q(g, g \wedge A) ,\ \ \ \ (b)\ \ Q(g, G) \ =\ 0 , 
\label{eqn2.1aa}
\end{eqnarray}
where 
$G = \frac{1}{2}\, g \wedge g$
and
$A$ is a symmetric $(0,2)$-tensor on $M$.
Using (\ref{eqn2.1}) and (\ref{eqn2.1aa})(a) we get
\begin{eqnarray}
Q(S,C) &=& Q\left( S, R - \frac{1}{n-2}\, g\wedge S + \frac{\kappa }{(n-2)(n-1)}\, G \right)\nonumber\\
&=& 
Q(S,R) - \frac{1}{n-2}\, Q(S, g\wedge S ) + \frac{\kappa }{(n-2)(n-1)}\,Q(S,G)\nonumber\\
&=& Q(S,R) + \frac{1}{n-2}\, Q\left( g, \frac{1}{2}\, S \wedge S \right) - \frac{\kappa }{(n-2)(n-1)}\,Q( g, g \wedge S) .
\label{eqn2.1srsr}
\end{eqnarray}
\begin{lem}  
(cf. {\cite[Lemma 1]{1989_DDV}})
(i) 
If $A$ is a symmetric $(0,2)$-tensor and $T$ a generalized curvature tensor
at a point $x$ of a semi-Riemannian manifold $(M,g)$, $n \geq 3$,
then the tensor $Q(A,T)$ satisfies at this point the identity
\begin{eqnarray*}
\sum _{ (X_{1}, X_{2}), (X_{3}, X_{4}), (X_{5},X_{6}) } Q(A,T)(X_{1},X_{2},X_{3},X_{4},X_{5},X_{6})
\ =\ 0 .
\end{eqnarray*}
(ii)
If $T$ is a generalized curvature tensor 
at a point $x$ of a semi-Riemannian manifold $(M,g)$, $n \geq 3$,
then the tensor $Q(g,T)$ vanishes at $x$
if and only if $T = \frac{\kappa(T)}{(n-1)n}\, G$ at this point, 
where $\kappa(T)$ is the scalar curvature of $T$.
\end{lem}
\begin{prop}
{\cite[Theorem 3.4 (i)] {DGJZ-2016}}
On any semi-Riemannian manifold $(M,g)$, $n \geq 4$, the following identity is satisfied 
\begin{eqnarray}
R \cdot C + C \cdot R 
&=& 
R \cdot R + C \cdot C - \frac{1}{(n-2)^{2}}\, Q\left( g,  - \frac{\kappa }{ n-1}\, g\wedge S + g \wedge S^{2} \right) .
\label{identity01}
\end{eqnarray}
\end{prop}
\begin{thm}
Let $(M,g)$, $n \geq 4$, be a semi-Riemannian Einstein manifold. We have
\newline
(i) {\cite[Section 5] {DGHJZ01}} The condition (\ref{Roterformula}) holds on $M$.
\newline
(ii) {\cite[Theorem 3.1] {DGJZ-2016}} If the condition (\ref{pseudo}) 
is satisfied on ${\mathcal{U}}_{R} \subset M$ then on this set we have
\begin{eqnarray*}
R \cdot R  &=& Q(S,R) + \left( L_{R} - \frac{\kappa}{n} \right) Q(g,C),\\  
C \cdot C  &=& \left( L_{R} - \frac{\kappa}{(n-1) n} \right) Q(g,C),\\
R \cdot C + C \cdot R &=& Q(S,C) + \left( 2 L_{R} - \frac{\kappa }{n-1} \right) Q(g,C) .
\end{eqnarray*}
\end{thm}
{\bf{Proof.}} (i) 
On any Einstein manifold $(M,g)$, $n \geq 4$, 
the following identity is satisfied (see, e.g.,  {\cite[p. 107] {DP-TVZ}})
$R \cdot C - C \cdot R = ( \kappa / ((n-1)n) ) \, Q(g,C)$.
This, by making use of 
$\kappa / ( (n-1)n ) = \kappa /(n-1) - \kappa /n$
and (\ref{einstein000}) 
turns into (\ref{Roterformula}), completing the proof. \qed
\begin{thm} 
Let $(M,g)$, $n \geq 4$, be a semi-Riemannian manifold satisfying (\ref{eq:h7a})
on ${\mathcal U}_{S} \cap {\mathcal U}_{C} \subset M$. We have
\newline
(i) {\cite{{DGHSaw}, {G5}}}
The following equations are satisfied on
${\mathcal U}_{S} \cap {\mathcal U}_{C}$:
\begin{eqnarray*}
S^{2} &=& \alpha _{1}\, S + \alpha _{2} \, g ,
\ \ \
\alpha _{1} 
\ =\ 
\kappa + ((n-2)\mu -1 ) \phi ^{-1} ,
\ \ \
\alpha _{2}
\ =\
( \mu \kappa + (n-1) \eta ) \phi ^{-1} ,\\ 
R \cdot C &=& L_{R}\, Q(g,C),
\ \ \
L_{R}
\ =\
 ((n-2) (\mu ^{2} - \phi \eta) - \mu ) \phi ^{-1} ,\\
R \cdot R &=& L_{R}\, Q(g,R),
\ \ \
R \cdot S \ =\ L_{R}\, Q(g,S),
\end{eqnarray*}
\begin{eqnarray*}
R \cdot R &=& Q(S,R) + L \, Q(g,C) ,
\ \ \
L
\ =\
L_{R} + \mu \phi ^{-1}
\ =\
(n-2) (\mu ^{2} - \phi \eta) \phi ^{-1},\\
C \cdot C &=& L_{C}\, Q(g,C) ,
\ \ \ 
L_{C} 
\ =\ L_{R} + \frac{1}{n-2} \left( \frac{\kappa }{n-1} - \alpha _{1} \right) ,\\
C \cdot R &=& L_{C}\, Q(g,R) , \ \ \  
C \cdot S \ =\ L_{C}\, Q(g,S) ,
\end{eqnarray*}
\begin{eqnarray*}
R \cdot C - C \cdot R &=& \frac{1}{n-2}\, Q(S,R)
+ \left( \frac{(n-1) \mu - 1}{(n-2) \phi} + \frac{\kappa }{n-1} \right) Q(g,R)\\
& &
+ \frac{\mu ((n-1) \mu - 1) - (n-1) \phi \eta }{(n-2) \phi} \, Q(S,G) ,
\end{eqnarray*}
\begin{eqnarray*}
R \cdot C - C \cdot R &=& 
\left( \frac{1}{\phi} \left( \mu - \frac{1}{n-2} \right) + \frac{\kappa }{n-1} \right) Q(g,R)
+ \left( \frac{\mu}{\phi } \left( \mu - \frac{1}{n-2}\right) - \eta \right) Q(S,G) .
\end{eqnarray*}
(ii) {\cite[Theorem 3.2 and Proposition 3.3] {DGJZ-2016}}
Moreover, we also have on ${\mathcal U}_{S} \cap {\mathcal U}_{C}$: 
(\ref{Roterformula}) and
\begin{eqnarray*}    
R \cdot C + C \cdot R  
&=& Q(S,C) + \left( L + L_{C} - \frac{1}{ (n-2) \phi } \right) Q(g,C) .
\end{eqnarray*}
\end{thm}

\noindent
{\bf{Remark 2.5.}} 
(i) 
(see {\cite{{Ch-DDGP}, {2016_DGHZhyper}, {DGJZ-2016}}} and references therein)
It is known that every warped product manifold $\overline{M} \times _{F} \widetilde{N}$
with an $1$-dimensional base manifold $(\overline{M}, \overline{g})$ and
a $2$-dimensional manifold $(\widetilde{N}, \widetilde{g})$
or an $(n-1)$-dimensional Einstein manifold
$(\widetilde{N}, \widetilde{g})$, $n \geq 4$, and a warping function $F$,
is a  quasi-Einstein manifold 
satisfyng (\ref{Riccipseudo07}).
\newline
(ii) 
It is easy to see that on the set ${\mathcal U}_{S}$ of any manifold $(M,g)$ the condition
(\ref{quasi02}) is equivalent  to
$(S - \alpha\, g) \wedge (S - \alpha\, g) = 0$.  
This gives
$(1/2)\, S \wedge S = \alpha\, g \wedge S - \alpha ^{2}\, G$.
From the last equation, by making use of (\ref{eqn2.1aa})(b), we obtain
$Q(g, (1/2)\, S \wedge S ) = \alpha\, Q(g, g \wedge S )$. 
This and (\ref{eqn2.1srsr}) yield
\begin{eqnarray}
Q(S,R) &=& Q(S,C) - \frac{1}{n-2}\, \left( \alpha - \frac{\kappa }{n-1} \right) Q( g, g \wedge S) .
\label{eqn2.1trtrtr}
\end{eqnarray}
(iii) 
(a)
(see {\cite[Example 4.1] {DGJZ-2016}} and references therein) 
The warped product manifold $\overline{M} \times _{F} \widetilde{N}$
with an $1$-dimensional base manifold $(\overline{M}, \overline{g})$,
$\overline{g}_{11} = \pm 1$,
and an $(n-1)$-dimensional Einstein fiber $(\widetilde{N}, \widetilde{g})$, $n \geq 5$, 
which is not a space of constant curvature, 
and a warping function $F$, 
satisfies on 
${\mathcal U}_{S} \cap {\mathcal U}_{C} \subset \overline{M} \times _{F} \widetilde{N}$
the following conditions: (\ref{Riccipseudo07}), with some function $L_{S}$, 
(\ref{quasi02}), with $\alpha = (\kappa / (n-1)) - L_{S}$, and (\ref{quasi022}).
In particular,
if $\overline{g}_{11} = - 1$ and $(\widetilde{N}, \widetilde{g})$, $n \geq 3$,
is a Riemannian manifold then
$\overline{M} \times _{F} \widetilde{N}$
is a special generalized Robertson-Walker spacetime \cite{ARS-1995}.
Generalized Robertson-Walker spacetimes satisfying 
curvature conditions of pseudosymmetry type
were investigated among others in
{\cite{{P119}, {Ch-DDGP}, {DGJZ-2016}, {DeScher}}}.
We also mention that 
Einstein generalized Robertson-Walker spacetimes 
were classified in \cite{ARS-1997}.
\newline
(b) 
(see \cite{AD-2002} and {\cite[Example 7.5 (i)] {2016_DGHZhyper}})
Let $M$ be a hypersurface in an Euclidean space 
${\mathbb{E}}^{n+1}$, $n = 2p + 1$, $p \geq 2$, 
having at every point three principal curvatures 
$\lambda _{1} = \lambda \neq 0$, $\lambda _{2} = - \lambda$ and $\lambda _{3} = 0$, 
provided that the multiplicity of $\lambda _{1}$ and $\lambda _{2}$ is $p$.
Evidently, we have ${\mathcal U}_{H} = M$. We can check that
the following conditions are satisfied on ${\mathcal U}_{H}$:
$H^{3} = - ( \kappa / (n-1) )\, H$, 
(\ref{Riccipseudo07}), with $L_{S} = 0$, 
(\ref{quasi02}), with $\alpha = \kappa / (n-1)$, 
$R \cdot C = Q(S,C)$, $C \cdot R = ( (n-3)/(n-2) ) Q(S,C)$,
and, 
$(n-2) (R \cdot C - C \cdot R) = Q(S,C)$,
i.e. (\ref{quasi022}), with $L_{S} = 0$.
\newline
(iv) In the standard Schwarzschild coordinates $(t; r; \theta; \phi)$, 
and the physical units ($c = G = 1$), the Reissner-Nordstr\"{o}m-de Sitter ($\Lambda > 0$), and 
the Reissner-Nordstr\"{o}m-anti-de Sitter ($\Lambda < 0$) metrics are given by the line element (see, e.g., \cite{SH})
\begin{eqnarray}
\ \ \ \ \ \ \ ds^{2} &=& - h(r)\, dt^{2} + h(r)^{-1}\, dr^{2} + r^{2}\, ( d\theta^{2} + \sin ^{2}\theta \, d\phi^{2}),\ \
h(r) \ =\ 1 - \frac{2M}{r} + \frac{Q^{2}}{r^{2}} - \frac{\Lambda r^{2}}{3}  ,
\label{rns01}
\end{eqnarray}
where $M$, $Q$ and $\Lambda$ are non-zero constants. 
\newline
(v) {\cite[Section 6] {DGHJZ01}} 
The metric (\ref{rns01}) satisfies (\ref{eq:h7a}) with 
\begin{eqnarray*}
\phi &=&
\frac{3}{2}\, ( Q^{2}  -  M r) r^{4} Q^{-4} ,\ \ \ \mu \ =\
\frac{1}{2}\, ( Q^{4} + 3 Q^{2} \Lambda r^{4} - 3 \Lambda M r^{5} ) Q^{-4} ,\\
\eta
&=&
\frac{1}{12}\, ( 3 Q^{6}
+ 4 Q^{4} \Lambda r^{4}
- 3 Q^{4} M r
+ 9 Q^{2} \Lambda ^{2} r^{8}
- 9 \Lambda^{2} M r^{9} ) r^{-4} Q^{-4}   .
\end{eqnarray*}
If we set $\Lambda = 0$ in (\ref{rns01})  
then we obtain the line element of the Reissner-Nordstr\"{o}m spacetime, 
see, e.g., {\cite[Section 9.2] {GrifPod}} and references therein. 
It seems that the Reissner-Nordstr\"{o}m spacetime is the "oldest" example 
of the Roter type warped product manifold.
\newline
(vi) 
Some comments on pseudosymmetric manifolds (also called Deszcz symmetric spaces),
as well as Roter spaces, are given in {\cite[Section 1] {DecuP-TSVer}}:
"{\sl{From a geometric point of view, the Deszcz symmetric spaces may well
be considered to be the simplest Riemannian manifolds next to the real space forms.}}" 
and 
"{\sl{From an algebraic point of view, Roter spaces may well be considered to
be the simplest Riemannian manifolds next to the real space forms.}}"
For further comments we refer to \cite{LV3-Foreword}.
\newline
(vii) 
Very recently, Roter spaces admitting geodesic mappings were studied in \cite{DH-2018}.

\section{Hypersufaces in space forms}

Let $M$, $n \geq 3$, be a connected hypersurface 
isometrically immersed in a semi-Riemannian space of constant curvature 
$N_{s}^{n+1}(c)$, with signature $(s,n+1-s)$. 
The Gauss equation of $M$ in $N_{s}^{n+1}(c)$ reads
\begin{eqnarray}
R &=& \frac{\varepsilon }{2}\, H \wedge H + \frac{\widetilde{\kappa}}{n(n+1)} \, G,\ \ \  
G \ =\ \frac{1}{2}\, g \wedge g ,\ \ \ 
\varepsilon \ =\ \pm 1 .
\label{C5}
\end{eqnarray}
From (\ref{C5}), by suitable contractions, we get
\begin{eqnarray}
S &=& \varepsilon \, ( \mathrm{tr}\, (H)\, H - H^{2} ) + \frac{ (n-1) \widetilde{\kappa}}{n(n+1)} \, g ,
\ \ \
\kappa \ =\ \varepsilon \, ( ( \mathrm{tr}\, (H))^{2} - \mathrm{tr}\, (H^{2}) ) 
+ \frac{ (n-1) \widetilde{\kappa}}{n+1} .
\label{C5ab}
\end{eqnarray}
It is known that (\ref{genpseudo01}) holds on $M$.
Precisely, we have on $M$
(see, e.g., {\cite[eq. (14)] {DGP-TV02}}, {\cite[Proposition 3.1] {DV-1991}})
\begin{eqnarray}
R \cdot R &=& Q(S,R) - \frac{(n-2) \widetilde{\kappa} }{n(n+1)}\, Q(g,C).
\label{900ab}
\end{eqnarray}
It is easy to see that (\ref{900ab}), 
by making use of (\ref{eqn2.1}) and (\ref{eqn2.1aa}), turns into
\begin{eqnarray}
R \cdot R  &=&
Q(S,R) - \frac{(n-2) \widetilde{\kappa} }{n(n+1)}\, Q(g,R) - \frac{ \widetilde{\kappa} }{n(n+1)}\, Q(S,G) . 
\label{900abdzdz}
\end{eqnarray}

We present now some results on pseudosymmetric hypersurfaces in $N_{s}^{n+1}(c)$, $n \geq 3$.

\begin{thm}
{\cite[Lemma 1, Theorem 1] {DVY}}
If $M$ is a hypersurface in $N_{s}^{n+1}(c)$, $n \geq 3$, 
satisfying (\ref{two}) on $\mathcal{U}_{R} \subset M$,
for some functions $\alpha $ and $\beta $, then on this set we have
\begin{eqnarray*} 
R \cdot R &=& \left( \frac{ \widetilde{\kappa} }{n(n+1)} - \varepsilon  \beta \right) Q(g,R) .
\end{eqnarray*}
\end{thm}
\begin{thm}
(i) {\cite[Remark 3.2] {D-1997}} 
If $M$ is a hypersurface in $N_{s}^{n+1}(c)$, $n \geq 4$, 
satisfying (\ref{pseudo}) on $\mathcal{U}_{R} \subset M$, 
then the following condition is satisfied on $\mathcal{U}_{R}$
\begin{eqnarray}
Q \left( S - \left( L_{R} + \frac{ (n-2) \widetilde{\kappa} }{n (n+1)} \right)  g , 
R - \frac{ \widetilde{\kappa} }{n (n+1) }\, G \right) &=& 0 .
\label{pseudo77}
\end{eqnarray}
(ii) 
{\cite[Lemma 1, Theorem 1] {D-1997}}
A hypersurface $M$ in 
$N_{s}^{n+1}(c)$, $n \geq 4$, is pseudosymmetric 
if and only if 
at every point of $M$ one of the following conditions is satisfied: 
(\ref{two}) or $\mathrm{rank}\, H = 2$.
\newline
(iii)
{\cite[Theorem 4.2] {1995_DDDVY}}
If $\mathrm{rank}\, H = 2$ at a point of a hypersurface $M$
in $N_{s}^{n+1}(c)$, $n \geq 4$, then (\ref{pseudo-constant01}) 
holds at this point.
\newline
(iv)
If $\mathrm{rank}\, H = 2$ at a point of ${\mathcal U}_{H} \subset M$
of a hypersurface $M$ 
in $N_{s}^{n+1}(c)$, $n \geq 4$, then at this point we have
\begin{eqnarray}
Q \left( S - \frac{ (n-1) \widetilde{\kappa} }{n (n+1)} \, g , R - \frac{ \widetilde{\kappa} }{n (n+1) }\, G \right) &=& 0 .
\label{pseudo88}
\end{eqnarray}
(v) (cf. {\cite[Section 4, eq. (62)] {DeScher}}, {\cite[Remark 5.2] {Kow02}}) 
If $\mathrm{rank}\, H = 2$ 
and $\mathrm{rank} \left( S - \frac{(n-1) \widetilde{ \kappa }}{n (n+1)}\, g \right) > 1$ 
at a point of ${\mathcal U}_{H} \subset M$ 
of a hypersurface $M$ 
in $N_{s}^{n+1}(c)$, $n \geq 4$, then at this point we have
\begin{eqnarray}
R - \frac{\widetilde{\kappa }}{n (n+1)}\, G
&=& \frac{\phi}{2} \left( S - \frac{ (n-1) \widetilde{ \kappa }}{n (n+1)}\, g \right)
\wedge \left( S - \frac{(n-1) \widetilde{ \kappa }}{n (n+1)}\, g \right)  ,\ \ \phi \in {\mathbb{R}} .
\label{pseudo-constant06}
\end{eqnarray}
\end{thm}
{\bf{Proof.}} 
(iv) The condition (\ref{pseudo88}) follows immediately from (iii) and (\ref{pseudo77}). 
(v) The condition (\ref{pseudo-constant06}) is an immediate consequence 
of (iv) and {\cite[Proposition 2.4] {DGHHY}} (see also {\cite[Lemma 3.1] {R102}}).
\qed


Note that from (\ref{C5ab}) it follows that (\ref{two}) and $\mathrm{tr}\, (H) - \alpha \neq 0$
are satisfied at every point of $({\mathcal U}_{S} \cap {\mathcal U}_{C}) \setminus {\mathcal U}_{H}$. 
We have
\begin{thm}
{\cite[Proposition 3.3] {2005-G}}
If $M$ is a hypersurface in $N_{s}^{n+1}(c)$, $n \geq 4$,
satisfying (\ref{two})
on $({\mathcal U}_{S} \cap {\mathcal U}_{C}) \setminus {\mathcal U}_{H}$,
where $\alpha $ and $\beta $ are some functions on this set,
then (\ref{eq:h7a}) holds on 
$({\mathcal U}_{S} \cap {\mathcal U}_{C}) \setminus {\mathcal U}_{H}$, 
where the functions $\phi, \mu, \eta $ are defined by
$\phi = \varepsilon (\mathrm{tr}(H) - \alpha )^{- 2}$ and
\begin{eqnarray*}
\mu \ =\ 
- \phi \left( \frac{(n-1)\widetilde{\kappa}}{n(n+1)} - \varepsilon \beta \right) ,\ \ \ \ 
\eta \ =\ 
\phi 
\left( \frac{(n-1)\widetilde{\kappa}}{n(n+1)} - \varepsilon \beta \right)^{2}  
+ \frac{\widetilde{\kappa} }{n(n+1)} .
\end{eqnarray*}
\end{thm}
\begin{thm}
If (\ref{two}) is satisfied at every point 
of a hypersurface in $N_{s}^{n+1}(c)$, $n \geq 4$, 
then (\ref{Roterformula}) holds on $M$.
\end{thm}
{\bf{Proof.}}
Evidently, (\ref{Roterformula}) is satisfied at all points of $M$ at which $C = 0$.
From Theorem 3.3 (i) it follows that our assertion is also true at all points of $M$ at which $S = (\kappa / n)\, g$. 
Finally, in view of Theorem 2.4 (ii) and Theorem 3.3, (\ref{Roterformula}) holds on 
${\mathcal U}_{S} \cap {\mathcal U}_{C} \subset M$. Our theorem is thus proved.
\qed


As an immediate consequence of the last result we have the following theorem.
\begin{thm}
If $M$ is a hypersurface in a Riemannian space of constant curvature $N^{n+1}(c)$, $n \geq 4$,
having at every point at most two distinct principal curvatures then
(\ref{Roterformula}) holds on $M$.
\end{thm}

\vspace{2mm}

\noindent
{\bf{Example 3.6.}} (cf. {\cite[Proposition 3.4, Example 3.2, Corollary 3.1] {2005-G}}) 
Let $M =  S^{p} ( \sqrt{ (p/n) } ) \times S^{n-p} ( \sqrt{ ( (n-p)/n ) } )$ 
be the Clifford torus in the $(n+1)$-dimensional unit sphere $S^{n+1}(1)$, $n \geq 4$.
Thus we have $c = \widetilde{\kappa} / (n(n+1)) = 1$. 
In addition, we assume that $n \neq 2p$ and $2 \leq p \leq n - 2$. 
Thus at every point of $M$ the tensors $S - (\kappa/n)\, g$ and $C$ are non-zero. 
Therefore ${\mathcal U}_{S} \cap {\mathcal U}_{C} = M$.
As it was shown in {\cite[Proposition 3.4] {2005-G}},
(\ref{eq:h7a}) is satisfied on $M$. Precisely, we have on $M$
\begin{eqnarray*}
R &=& \frac{ p (n-p) }{ 2 (n - 2p)^{2} } \, ( S - (n-2)\, g) \wedge ( S - (n-2)\, g) + \frac{1}{2}\, g \wedge g . 
\end{eqnarray*}
This, in view of Theorem 2.2 (ii), implies (\ref{Roterformula}). Further,
it is known that every Clifford torus is a semisymmetric manifold. Thus $R \cdot R = 0$ and, in a consequence, 
$R \cdot C = 0$ on $M$.
Now (\ref{Roterformula}) reduces to
$C \cdot R = Q(S,C) - (\kappa / (n-1))\, Q(g,C)$ on $M$. 

\vspace{2mm}

\noindent
{\bf{Remark 3.7.}} Let $M$ be a hypersurface in $N_{s}^{n+1}(c)$, $n \geq 4$. 
(i) 
If 
(\ref{einstein000}) holds at a point $x \in M$ then at this point (\ref{C5ab}) turns into
\begin{eqnarray*}
H^{2} &=& \mathrm{tr}\, (H)\, H  
+ \frac{(n-1) \varepsilon \,}{n} \left( \frac{ \widetilde{\kappa}}{n+1}  - \frac{\kappa }{n-1} \right) g .
\end{eqnarray*}
(ii) 
The Weyl conformal curvature tensor $C$ of $M$ vanishes at a point $x \in M$ if and only if at $x$ we have
$\mathrm{rank}\, ( H - \alpha_{1}\, g) \leq 1$, for some $\alpha_{1} \in {\mathbb{R}}$ {\cite[Theorem 4.1] {DV-1991}}. 
But the last condition,
in view of {\cite[Lemma 2.2] {ADEMO-2002}},
implies (\ref{C5ab}), for some $\alpha , \beta \in {\mathbb{R}}$. 

\section{The condition $H^{3} = \mathrm{tr} (H)\, H^{2} + \psi \, H + \rho \, g$}

In this section we consider hypersurfaces $M$
in a semi-Rieman\-nian space of constant curvature $N_{s}^{n+1}(c)$, $n \geq 4$,
satisfying 
on ${\mathcal U}_{H} \subset M$
curvature conditions of the kind:
the tensor $R \cdot C$, $C \cdot R$ or $R \cdot C - R \cdot C$ is a linear combination
of the tensor $R \cdot R$ and of a finite sum of the Tachibana tensors of the form $Q(A,B)$, 
where $A$ is a symmetric $(0,2)$-tensor and $B$ a generalized curvature tensor.
As it was mentioned in Section 1, 
if such condition is satisfied on ${\mathcal U}_{H}$ then (\ref{DS4}) holds on this set.

\begin{prop} {\cite[Proposition 5.1, eq. (29)]{Saw-2005}}
If $M$ is a hypersurface 
in $N_{s}^{n+1}(c)$, $n \geq 4$,
satisfying (\ref{DS4}) on ${\mathcal U}_{H} \subset M$,
for some functions 
$\psi $ and $\rho $ on ${\mathcal U}_{H}$,
then on this set we have
\begin{eqnarray}
R \cdot C 
&=& 
Q(S,R) 
- \frac{(n-2) \widetilde{\kappa }}{n (n+1)}\, Q(g,R)
+ \alpha _{2}\, Q(S,G)
+ \frac{\rho }{n-2}\, Q(H,G) ,
\label{ZZ1}
\end{eqnarray}
\begin{eqnarray}
C \cdot R &=&
\frac{n-3}{n-2}\, Q(S,R) 
+ \alpha _{1}\, Q(g,R)
+ \alpha _{2}\,  Q(S,G) ,
\label{ZZ2}
\end{eqnarray}
\begin{eqnarray}
(n-2)\, (R \cdot C - C \cdot R) &=& Q(S,R) + \rho \, Q(H,G)\nonumber\\
& & + \left( \frac{(n-1) \widetilde{\kappa} }{n(n+1)} - \frac{\kappa }{n-1} 
- \varepsilon \psi \right) Q(g, R),
\label{ZZ3}
\end{eqnarray}
\begin{eqnarray}
(n-2)\, C \cdot C &=& 
(n-3)\, Q(S,R) + (n-2) \alpha _{1}\, Q(g,R)\nonumber\\
& &
+ ( \alpha _{1} - \alpha _{2} )\, Q(S,G)
+ \frac{n-3}{n-2} \rho \, Q(H,G) ,
\label{DS16A}
\end{eqnarray}
\begin{eqnarray}
R \cdot S &=& 
\frac{\widetilde{\kappa }}{n(n+1)}\, Q(g,S) + \rho \, Q(g,H) ,
\label{DZ004}
\end{eqnarray}
\begin{eqnarray}
\alpha _{1} 
&=&
\frac{1}{n-2} \left( \frac{\kappa }{n-1} + \varepsilon \psi
- \frac{(n^{2} - 3n + 3) \widetilde{\kappa }}{n(n+1)} \right) ,
\label{AA00two}\\
\alpha _{2} 
&=&
- \frac{(n-3) \widetilde{\kappa} }{(n-2)n (n+1)} .
\label{AA00}
\end{eqnarray}
\end{prop}
\begin{prop}
Let $M$ be a hypersurface in 
$N_{s}^{n+1}(c)$, $n \geq 4$,
satisfying (\ref{DS4}) on ${\mathcal U}_{H} \subset M$,
for some functions 
$\psi $ and $\rho $ on ${\mathcal U}_{H}$. We have
\newline
(i) The following conditions are satisfied on ${\mathcal U}_{H}$
\begin{eqnarray}
Q(\rho\, H - \alpha_{3}\, S - S^{2},G) &=& 0 ,
\label{DZ005}
\end{eqnarray}
\begin{eqnarray}
\rho\, H  &=& S^{2} + \alpha_{3}\, S  + \frac{\lambda}{n}\, g ,\ \ \ 
\lambda \ =\ \rho\, \mathrm{tr}(H) - \kappa \, \alpha_{3} - \mathrm{tr}(S^{2}) , 
\label{DZ008}
\end{eqnarray}
\begin{eqnarray}
\ \ \ \ \ \
\alpha _{3}
&=&
(n-2)^{2} 
\left( \frac{1}{n-2} (\alpha_{1} - \alpha_{2}) - 2\alpha_{2} - \frac{ \widetilde{\kappa} }{n(n+1)} \right) 
- \frac{\kappa }{n-1} \ =\
\varepsilon \psi - \frac{ 2 ( n-1) \widetilde{\kappa} }{n(n+1)} , 
\label{DZ006}
\end{eqnarray}
where 
$\alpha _{1}$ and $\alpha _{2}$ are defined by (\ref{AA00two}) and (\ref{AA00}), respectively.
Moreover, 
(\ref{GGG01}), (\ref{EEE01}) and (\ref{EEE01new}) hold on ${\mathcal U}_{H}$, where $\rho _{1}$, $\rho _{2}$ and $\rho _{3}$
are defined by
\begin{eqnarray}
\rho _{1} &=& - \frac{ (n-2) \widetilde{\kappa} }{n (n+1) } - \alpha _{3} ,\ \ \ \
\rho _{2} \ =\ - \frac{\lambda}{n} - \left( \frac{ (n-1) \widetilde{\kappa} }{n (n+1) } + \alpha _{3}\right)  \alpha _{3} ,\nonumber\\
\rho _{3} &=&
\frac{1}{n} 
\left( \mathrm{tr}(S^{3}) + \left( 2 \varepsilon \psi -    \frac{ 3 (n-1) \widetilde{\kappa} }{n (n+1) }      \right)
\mathrm{tr}(S^{2}) - \kappa \rho _{2} \right) ,
\label{EEE01bb}
\end{eqnarray}
respectively. 
\newline
(ii) If at a point $x \in {\mathcal U}_{H}$ we have
$S^{2} = \beta_{1}\, S + \beta _{2}\, g$, for some $\beta_{1}, \, \beta_{2} \in {\mathbb{R}}$, 
then $\rho = 0$, $\beta _{1} = \alpha _{3}$ and $\beta _{2} = - ( \lambda / n)$ at this point. 
\end{prop}
{\bf{Proof.}} 
(i) The identities (\ref{identity01}) and (\ref{900abdzdz}), 
together with (\ref{ZZ1}), (\ref{ZZ2}) and (\ref{DS16A}), give
\begin{eqnarray*}
& &
\left( 1 + \frac{n-3}{n-2} \right) Q(S,R) 
+ \alpha_{1}\, Q(g,R) + 2\alpha_{2}\, Q(S,G) 
- \frac{(n-2) \widetilde{\kappa}  }{n(n+1)}\, Q(g,R)
+ \frac{\rho }{n-2}\, Q(H,G)\\
&=&
\left( 1 + \frac{n-3}{n-2} \right) Q(S,R)
+ \alpha_{1}\, Q(g,R) 
+ \frac{1}{n-2}\left( (\alpha_{1} - \alpha_{2})\, Q(S,G) 
+  \frac{(n-3) \rho }{n-2}\, Q(H,G)\right)\\
& &
- \frac{(n-2) \widetilde{\kappa} }{n(n+1)}\, Q(g,R) 
- \frac{ \widetilde{\kappa} }{n(n+1)}\, Q(S,G)
-  \frac{1}{(n-2)^{2}} Q\left( g,  - \frac{\kappa }{ n-1}\, g\wedge S + g \wedge S^{2} \right) ,
\end{eqnarray*}
\begin{eqnarray*}
\frac{ \rho }{(n-2)^{2}}\, Q(H,G)
&=&
\left( \frac{1}{n-2} (\alpha_{1} - \alpha_{2}) - 2\alpha_{2} - \frac{ \widetilde{\kappa} }{n(n+1)}
 \right) Q(S,G)\\ 
& & -  \frac{1}{(n-2)^{2}}\, Q\left( g,  - \frac{\kappa }{ n-1}\, g\wedge S + g \wedge S^{2} \right) ,
\end{eqnarray*}
\begin{eqnarray*}
 \rho\, Q(H,G)
&=&
(n-2)^{2} 
\left( \frac{1}{n-2} (\alpha_{1} - \alpha_{2}) - 2\alpha_{2} - \frac{ \widetilde{\kappa} }{n(n+1)}  \right) 
Q(S,G)\\ 
& & -  Q\left( g,  - \frac{\kappa }{ n-1}\, g\wedge S + g \wedge S^{2} \right) .
\end{eqnarray*}
This, by making use of (\ref{eqn2.1aa}) and (\ref{DZ006}), yields 
$\rho\, Q(H,G)
= \alpha_{3}\, Q(S,G) -  Q(g, g \wedge S^{2})$ 
and
$\rho\, Q(H,G) 
= \alpha_{3}\, Q(S,G) + Q(S^{2}, G)$
and (\ref{DZ005}).
From (\ref{DZ005}), by a suitable contraction, we get
$Q(\rho\, H - \alpha_{3}\, S - S^{2},g) = 0$
and, in a consequence, (\ref{DZ008}). 
Now we prove that 
(\ref{GGG01}) and (\ref{EEE01}) hold on ${\mathcal U}_{H}$.
From (\ref{DZ008}) we get
\begin{eqnarray}
\rho \, Q(g,H) &=& Q(g, S^{2}) + \alpha _{3}\, Q(g,S) ,
\label{EEE02}\\
\rho \, Q(H,S) &=& - Q(S, S^{2}) + \frac{\lambda }{n}\, Q(g,S) ,
\label{EEE03}\\
R \cdot S^{2} &=& \rho \, (R \cdot H) - \alpha _{3}\, (R \cdot S) .
\label{EHEHEH03}
\end{eqnarray}
The conditions 
(\ref{DZ004}), (\ref{DZ006}) and (\ref{EEE02}) yield immediately (\ref{GGG01}).
Further,
(\ref{EHEHEH03}), by (\ref{DZ004}) and (\ref{EEE02}), turns into
\begin{eqnarray}
R \cdot S^{2} &=& \rho \, (R \cdot H)
-
\left(  \frac{ \widetilde{\kappa} \alpha _{3} }{n (n+1) } + \alpha _{3}^{2} \right) Q(g, S)
- \alpha _{3} \, Q(g, S^{2}) .
\label{EEE05}
\end{eqnarray}
From the Gauss equation (\ref{C5}) of $M$ in $N_{s}^{n+1}(c)$ we get
\begin{eqnarray*}
g^{rs}( H_{hr}R_{sijk} + H_{ir}R_{shjk} ) &=& 
\varepsilon \, 
( H^{2}_{hk}H_{ij} - H^{2}_{hj}H_{ik} + H^{2}_{ik}H_{hj} - H^{2}_{ij}H_{hk}) \nonumber\\
& &
+ \frac{ \widetilde{\kappa} }{n (n+1) }\,
( H_{hk}g_{ij} - H_{hj}g_{ik} + H_{ik}g_{hj} - H_{ij}g_{hk}) ,\\
R \cdot H &=& \varepsilon \, Q(H, H^{2}) + \frac{ \widetilde{\kappa}}{ n (n+1) } \, Q(g,H).
\end{eqnarray*}
This, together with 
\begin{eqnarray*}
\varepsilon \, Q(H, H^{2}) &=& Q(H,S) - \frac{ \widetilde{\kappa}}{ n (n+1)} \, Q(g,H) , 
\end{eqnarray*}
which 
is an immediate consequence of (\ref{C5ab}), turns into
\begin{eqnarray*}
R \cdot H &=&  - Q(H,S) - \frac{ (n - 2) \widetilde{\kappa}}{ n (n+1)} \, Q(g,H) .
\end{eqnarray*}
Applying in this (\ref{EEE02}) and (\ref{EEE03}) we obtain 
\begin{eqnarray*}
\rho \, ( R \cdot H ) &=& Q(S, S^{2}) 
- 
\frac{ (n-2) \widetilde{\kappa} }{n (n+1) } \, Q(g, S^{2})
-
\left( \frac{\lambda}{n} + \frac{ (n-2) \widetilde{\kappa} \alpha _{3} }{n (n+1) }  \right) Q(g, S ) ,
\end{eqnarray*}
which together with (\ref{EEE05}) yields (\ref{EEE01}).
We prove now that (\ref{EEE01new}) holds on ${\mathcal U}_{H}$.
From (\ref{GGG01}) it follows that
\begin{eqnarray}
S_{h}^{r}R_{rijk} + S_{i}^{r}R_{rhjk} &=&
\rho_{4}\, ( g_{hj}S_{ik} + g_{ij}S_{hk} - g_{hk}S_{ij} - g_{ik}S_{hj} )\nonumber\\
& &
+ g_{hj}S^{2}_{ik} + g_{ij}S^{2}_{hk} - g_{hk}S^{2}_{ij} - g_{ik}S^{2}_{hj} ,
\label{1240}
\end{eqnarray}
where
$S_{h}^{r} = S_{hk}g^{kr}$
and 
$\rho _{4} = \varepsilon \psi - (( 2 n - 3) \widetilde{\kappa })/ (n(n+1))$.
Transvecting (\ref{1240}) with $S_{l}^{h}$ we get
\begin{eqnarray*}
S_{l}^{h} S_{h}^{r} \, R_{rijk} + S_{l}^{h} S_{i}^{r}R_{rhjk} &=&
\rho_{4}\, ( S_{lj}S_{ik} - S_{lk}S_{ij} + g_{ij}S^{2}_{lk}  - g_{ik}S^{2}_{lj} )\\
& &
+ S_{lj}S^{2}_{ik} - S_{lk}S^{2}_{ij} + g_{ij}S^{3}_{hk}  - g_{ik}S^{3}_{lj} ,
\end{eqnarray*}
which by symmetrization in $\, i,l\, $ leads to
$R \cdot S^{2} = Q(S, S^{2}) + Q(g, S^{3}) +  \rho_{4}\, Q(g,S^{2})$.
This and (\ref{EEE01}) yield 
$Q(g, S^{3} + (\rho _{4} - \rho _{1})\, S^{2} - \rho _{2}\, S ) = 0$,
which in view of {\cite[Lemma 2.4 (i)] {DV-1991}} implies 
(\ref{EEE01new}),  completing the proof of (i).
(ii) From (\ref{DZ008}), by our assumption and (\ref{C5ab}), we obtain
\begin{eqnarray*}
\rho\, H  &=& (\beta _{1} - \alpha_{3})\, \varepsilon \, ( \mathrm{tr}\, (H)\, H - H^{2} )
+
\left( (\beta _{1} - \alpha_{3})\,  \frac{ (n-1) \widetilde{\kappa}}{n(n+1)} +  \beta _{2} + \frac{\lambda}{n} \right) g ,
\end{eqnarray*}
which completes the proof. 
\qed

\begin{prop}
Let $M$ be a hypersurface in 
$N_{s}^{n+1}(c)$, $n \geq 4$. We have
\newline
(i) {\cite[Proposition 2.1] {2005-DSaw}} If $n = 4$ then
(\ref{DS4}) reduces on ${\mathcal U}_{H} \subset M$ to (\ref{DS4aa}).
\newline
(ii)
{\cite[Theorem 3.2 (iv), Proposition 3.1] {Saw-2004}}
If $\mathrm{rank}\, H = 2$
at every point of 
${\mathcal U}_{H} \subset M$ then (\ref{DS4aa}) is satisfied on this set with 
\begin{eqnarray}
\psi 
&=&  
\frac{1}{2}\, (  \mathrm{tr}\, (H^{2}) - ( \mathrm{tr}\, (H) )^{2} )
\ =\ 
\frac{(n-1) \varepsilon }{2} \left( \frac{\widetilde{\kappa }}{n+1} - \frac{\kappa }{n-1} \right) .
\label{pseudo-constant05}
\end{eqnarray}
(iii) 
If at every point of 
${\mathcal U}_{H} \subset M$ the conditions:
$\mathrm{rank}\, H = 2$ and $\mathrm{rank} \left( S - \frac{(n-1) \widetilde{ \kappa }}{n (n+1)}\, g \right) > 1$
are satisfied then 
(\ref{eq:h7a}) holds 
on this set with the functions $\phi$, $\mu$ and $\eta$ defined by
\begin{eqnarray}
\frac{2}{(n-1) \phi } &=&  \frac{\widetilde{\kappa }}{n+1} - \frac{\kappa }{n-1} ,
\label{pseudo-constant14}\\
\mu &=&  - \frac{ (n-1) \widetilde{\kappa } }{ n (n+1) }\, \phi ,\ \
\eta \ =\  \frac{ \widetilde{\kappa } }{ n (n+1) } 
\left(
\frac{ (n-1)^{2} \widetilde{\kappa } }{ n (n+1) }\, \phi  + 1 
\right) ,
\label{pseudo-constant15}
\end{eqnarray}
respectively. Moreover, 
the following conditions are satisfied on ${\mathcal U}_{H}$: 
(\ref{pseudo-constant01}) and
\begin{eqnarray}
R \cdot C &=&  \frac{\widetilde{\kappa }}{n (n+1)}\, Q(g,C) ,
\label{pseudo-constant02}\\
C \cdot R &=&
\frac{n-3}{(n-2) (n-1) \phi } \, Q(g,R) ,
\label{pseudo-constant03}\\
C \cdot C &=&
\frac{n-3}{(n-2) (n-1) \phi } \, Q(g,C) ,
\label{pseudo-constant04}
\end{eqnarray}
(cf., {\cite[Theorem 3.2(ii), Proposition 4.3] {Saw-2004}}). 
\newline
(iv) (cf., {\cite[Theorem 3.1] {ADY-1997}}) 
On ${\mathcal U}_{H} \subset M$ 
(\ref{pseudo-constant01}) and
(\ref{pseudo-constant02}) are equivalent.
\newline
(v) {\cite[Proposition 3.2] {1995_DDDVY}}
If (\ref{DS4aa}) is satisfied on  ${\mathcal U}_{H} \subset M$ then (\ref{DZ004Ricciaaa})
holds on this set.
\end{prop}
{\bf{Proof.}} 
(iii) From Theorem 3.2 (v) it follows that (\ref{pseudo-constant15}) holds on ${\mathcal U}_{H}$, 
where $\phi $ is some function on this set. We prove now that $\phi $ satisfies (\ref{pseudo-constant14}).
From Theorem 2.4 (i) and (\ref{pseudo-constant15}) it follows that
\begin{eqnarray}
S^{2} &=& 
\left( \kappa - \frac{ (n-2) (n-1) \widetilde{\kappa } }{n (n+1)} - \phi ^{-1} \right) S
- \frac{ (n-1) \widetilde{\kappa } }{n (n+1)} 
\left( \kappa - \frac{ (n-1)^{2} \widetilde{\kappa } }{n (n+1)} - \phi ^{-1} \right) g.
\label{pseudo-constant16}
\end{eqnarray}
But on the other hand, (\ref{DZ008}), by (\ref{DZ006}) making use of and (\ref{pseudo-constant05}), yields
\begin{eqnarray}
S^{2} &=& \frac{1}{2} 
\left( \kappa - \frac{ (n-1) (n-4) \widetilde{\kappa } }{n (n+1)} \right) S 
- \frac{\lambda }{n}\, g .
\label{pseudo-constant17}
\end{eqnarray}
Now using (\ref{pseudo-constant16}) and (\ref{pseudo-constant17})
we get (\ref{pseudo-constant14}). Further, using 
(\ref{pseudo-constant14}), (\ref{pseudo-constant15}) and suitable formulas given in Theorem 2.4 (i)
we can check that (\ref{pseudo-constant01}) 
and (\ref{pseudo-constant02})-(\ref{pseudo-constant04}) hold good, 
which completes the proof of (iii). 
\qed

\begin{thm}{\cite[Theorem 5.1, Theorem 5.2]{2011-DGPSS}}
Let $M$ be a hypersurface in 
$N_{s}^{n+1}(c)$, $n \geq 4$. We have
\newline
(i) If a generalized curvature tensor $B_{1}$ satisfies 
\begin{eqnarray}
R \cdot R &=& Q(g,B_{1})
\label{dgpss777}
\end{eqnarray}
on ${\mathcal{U}}_{H} \subset M$, then on this set we have
\begin{eqnarray}
(n-1)\, B_{1} &=& 
\left( \kappa + \varepsilon \psi  - \frac{ (n-1)^{2} \widetilde{\kappa} }{n(n+1)} \right) R 
- \frac{1}{2}\, S \wedge S + g \wedge S^{2}
\nonumber\\
& &
+ \left( \varepsilon \psi - \frac{(n - 1)\widetilde{\kappa}}{n(n+1)} \right) g \wedge S + \lambda \, G  ,
\label{dgpss26}
\end{eqnarray}
where $\lambda $ is some function on ${\mathcal{U}}_{H}$.
\newline
(ii) If a generalized curvature tensor $B_{2}$ satisfies 
\begin{eqnarray}
R \cdot C &=& Q(g,B_{2})
\label{dgpss201}
\end{eqnarray}
on ${\mathcal{U}}_{H} \subset M$, then on this set we have
\begin{eqnarray}
(n-1)\, B_{2} &=& 
\left( \kappa + \varepsilon \psi  - \frac{ (n-1)^{2} \widetilde{\kappa} }{n(n+1)}\right)
R 
-
\frac{1}{n-2}\, g \wedge S^{2} 
\nonumber\\
& &
-
\frac{1}{2}\, S \wedge S 
- \frac{1}{n-2}
\left( \varepsilon \psi - \frac{ (n-1)^{2} \widetilde{ \kappa} }{n(n+1)} \right)  g \wedge S
+ \lambda \, G  ,
\label{dgpss202}
\end{eqnarray}
where $\lambda $ is some function on ${\mathcal{U}}_{H}$.
\newline
(iii) 
If a generalized curvature tensor $B_{3}$ satisfies 
\begin{eqnarray}
C \cdot R &=& Q(g,B_{3})
\label{dgpss201a}
\end{eqnarray}
on ${\mathcal{U}}_{H} \subset M$, then on this set we have
\begin{eqnarray}
B_{3} &=& 
\left( \frac{ \kappa }{n-1} + \frac{ 2 \varepsilon \psi}{n-1} - \frac{  \widetilde{\kappa} }{n+1} \right) R 
+ \lambda \, G  \nonumber\\ 
& &
+ \frac{n-3}{(n-2)(n-1)} 
\left( \left( \varepsilon \psi - \frac{ (n-1) \widetilde{\kappa} }{n(n+1)} \right) g \wedge S 
- \frac{1}{2}\, S \wedge S + g \wedge S^{2} \right) ,
\label{dgpss202a}
\end{eqnarray}
where $\lambda $ is a function on ${\mathcal{U}}_{H}$.
(iv) If a generalized curvature tensor $B_{4}$ satisfies  
\begin{eqnarray}
R \cdot C - C \cdot R &=& Q(g,B_{4})
\label{dgpss201b}
\end{eqnarray}
on ${\mathcal{U}}_{H} \subset M$, then on this set we have
\begin{eqnarray}
B_{4} 
&=& 
\left( - \frac{\varepsilon \psi  }{n-1} + \frac{ \widetilde{\kappa} }{n (n+1)} \right) R
+
\left( - \frac{\varepsilon \psi  }{n-1} + \frac{ 2  \widetilde{\kappa} }{n (n+1)} \right)
g \wedge S
\nonumber\\
& &
- \frac{1}{n-1} \, g \wedge S^{2}
- \frac{1}{2 (n-2) (n-1) } \, S \wedge S + \lambda\, G ,
\label{dgpss202b}
\end{eqnarray}
where $\lambda $ is a function on ${\mathcal{U}}_{H}$.
\end{thm}
\begin{thm}
Let $M$ be a hypersurface in 
$N_{s}^{n+1}(c)$, $n \geq 4$,
satisfying (\ref{DS4}) on ${\mathcal{U}}_{H} \subset M$,
for some functions 
$\psi $ and $\rho $ on ${\mathcal{U}}_{H}$.
If the tensor $C \cdot C$ and a generalized curvature tensor $B$ satisfy 
\begin{eqnarray}
C \cdot C &=& Q(g,B)
\label{dgpss777cc}
\end{eqnarray}
on ${\mathcal{U}}_{H}$, then on this set we have
\begin{eqnarray}
B &=&
\left( \frac{ \kappa }{n-1} + \frac{ 2 \varepsilon \psi}{n-1} - \frac{  \widetilde{\kappa} }{n+1} \right) C  + \lambda \, G\nonumber\\
& &
- \frac{n-3}{(n-2)^{2} (n-1)} \left( \frac{n-2}{2}\, S \wedge S - \kappa\, g \wedge S + g \wedge S^{2} \right) ,
\label{dgpss777ccdd}
\end{eqnarray}
where $\lambda$ is some function on ${\mathcal{U}}_{H}$.
\end{thm}
{\bf{Proof.}} 
From (\ref{eqn2.1aa})(a), (\ref{DS16A}) and (\ref{dgpss777cc}) 
it follows that the tensor $Q(S,R)$ is a linear combination
of the tensors of the form $Q(g, g \wedge A)$, 
where $A$ is a symmetric $(0,2)$-tensor, and the tensors $Q(g,R)$ and $Q(g,B)$.
Now, from (\ref{900abdzdz}), (\ref{ZZ1}) and (\ref{ZZ2}) it follows 
that the tensors $R \cdot R$, $R \cdot C$ and  $C \cdot R$
satisfy (\ref{dgpss777}), (\ref{dgpss201}) and (\ref{dgpss201a}), respectively, 
where $B_{1}$, $B_{2}$ and  $B_{3}$ are some generalized curvature tensors.
In view of Theorem 4.4, 
the tensors $B_{1}$, $B_{2}$ and  $B_{3}$ satisfy 
(\ref{dgpss26}), (\ref{dgpss202}) and (\ref{dgpss202a}), respectively.
Further,
(\ref{dgpss777}), (\ref{dgpss201}) and (\ref{dgpss201a})
together with 
(\ref{identity01}) yield
\begin{eqnarray*}
C \cdot C 
&=&
\frac{1}{(n-2)^{2}}\, Q \left( g,  - \frac{\kappa }{ n-1}\, g\wedge S + g \wedge S^{2}\right)
+ Q(g, - B_{1} + B_{2} + B_{3}) .
\end{eqnarray*}
This by making use of (\ref{dgpss26}), (\ref{dgpss202}) and (\ref{dgpss202a}) 
leads to (\ref{dgpss777cc}), with
\begin{eqnarray}
& &
B \ =\ 
\left( \frac{ \kappa }{n-1} + \frac{ 2 \varepsilon \psi}{n-1} - \frac{  \widetilde{\kappa} }{n+1} \right) R 
- \frac{n-3}{2 (n-2) (n-1) } \, S \wedge S 
\nonumber\\ 
& &
- \frac{1}{n-2}\left( \frac{ \kappa }{(n-2)(n-1)} 
+ \frac{ 2 \varepsilon \psi}{n-1} - \frac{  \widetilde{\kappa} }{n+1} \right) g \wedge S
- \frac{n - 3}{(n-2)^{2} (n-1)} \, g \wedge S^{2} + \lambda _{2}\, G ,
\label{dgpss777ccddzz}
\end{eqnarray}
where $\lambda _{2}$ is some function on ${\mathcal{U}}_{H}$.
Now (\ref{dgpss777ccddzz}), by (\ref{eqn2.1}), gives (\ref{dgpss777ccdd}), 
completing the proof.
\qed


From Proposition 4.1 and theorems 4.3 and 4.4 it follows the following result.

\begin{thm}
Let $M$ be a hypersurface in 
$N_{s}^{n+1}(c)$, $n \geq 4$,
satisfying (\ref{DS4}) on ${\mathcal U}_{H} \subset M$,
for some functions 
$\psi $ and $\rho $ on ${\mathcal U}_{H}$.
If the tensor $Q(S,R)$ is equal to the Tachibana tensor $Q(g,T)$,
where $T$ is a generalized curvature tensor, then the tensors 
$R \cdot R$, $R \cdot C$, $C \cdot R$, $R \cdot C - C \cdot R$ and $C \cdot C$ satisfy
(\ref{dgpss777}),
(\ref{dgpss201}),
(\ref{dgpss201a}),
(\ref{dgpss201b}) and
(\ref{dgpss777cc}), with the tensors
$B_{1}$, $B_{2}$, $B_{3}$, $B_{4}$ and $B$, 
defined by
(\ref{dgpss26}),
(\ref{dgpss202}),
(\ref{dgpss202a}),
(\ref{dgpss202b})
and
(\ref{dgpss777ccdd}),
respectively.
\end{thm}
\begin{prop}
Let $M$ 
be a hypersurface in 
$N_{s}^{n+1}(c)$, $n \geq 4$. We have
\newline
(i) {\cite[Theorem 3.7]{DGJZ-2016}}
The following identity is satisfied on $M$
\begin{eqnarray}
R \cdot C + C \cdot R  
&=& Q(S,C) 
- \frac{(n-2) \widetilde{\kappa} }{n(n+1)}\, Q(g,C) 
+ C \cdot C \nonumber\\
& & - \frac{1}{(n-2)^{2}}\, Q \left( g, \frac{n-2}{2}\, S \wedge S - \kappa \, g\wedge S + g \wedge S^{2} \right) .
\label{02identity01hyper}
\end{eqnarray}
(ii) If (\ref{DS4}) is satisfied on ${\mathcal U}_{H} \subset M$,
for some functions 
$\psi $ and $\rho $,
then 
(\ref{DS16Anew01}) and (\ref{02identity01hyper17})
hold on this set.
\end{prop}
{\bf{Proof.}} 
(ii) From (\ref{ZZ1}), (\ref{DS16A}) and (\ref{AA00}) we get
\begin{eqnarray*}
& &
C \cdot C - \frac{n-3}{n-2}\, R \cdot C\nonumber\\
&=& 
\frac{n-3}{n-2}\, Q(S,R) + \alpha _{1}\, Q(g,R)
+ \frac{1}{n-2} \left( ( \alpha _{1} - \alpha _{2} )\, Q(S,G)
+ \frac{n-3}{n-2} \rho \, Q(H,G) \right)\nonumber\\
& &
- \frac{n-3}{n-2}\, Q(S,R) 
+ \frac{n-3}{n-2}\, \frac{(n-2) \widetilde{\kappa }}{n (n+1)}\, Q(g,R)
- \frac{n-3}{n-2}\, \alpha _{2}\, Q(S,G)
- \frac{n-3}{n-2}\, \frac{\rho }{n-2}\, Q(H,G)\nonumber\\
&=&
\left( \alpha_{1} + \frac{(n-3) \widetilde{\kappa }}{n (n+1)} \right) Q(g,R)
+ \left( \frac{\alpha _{1} - \alpha _{2}}{n-2} - \frac{n-3}{n-2}\alpha_{2} \right) Q(S,G)\nonumber\\
&=&
\left( \alpha_{1} + \frac{(n-3) \widetilde{\kappa }}{n (n+1)} \right) Q(g,R)
+ \left( \frac{\alpha _{1}}{n-2} - \frac{n - 3 + 1}{n-2}\alpha_{2} \right) Q(S,G)\nonumber\\
&=&
\alpha_{1}\, Q(g,R) + \frac{\alpha _{1}}{n-2}\, Q(S,G)
+ \frac{ (n-3) \widetilde{\kappa } }{n (n+1)} \, Q(g,R) 
- \alpha_{2}\, Q(S,G) 
\end{eqnarray*}
\begin{eqnarray*}
&=&
\alpha_{1}\, Q(g,C) + \frac{ (n-3) \widetilde{\kappa } }{n (n+1)} \, Q(g,R) 
+ \frac{ (n-3) \widetilde{\kappa } }{ (n-2) n (n+1)} \, Q(S,G) \nonumber\\
&=&
\alpha_{1}\, Q(g,C) + \frac{ (n-3) \widetilde{\kappa } }{n (n+1)}  
\left( Q(g,R) - \frac{1}{n-2}\, Q(g, g \wedge S) \right) \ =\
\left( \alpha_{1} + \frac{ (n-3) \widetilde{\kappa } }{n (n+1)} \right) Q(g,C)
\end{eqnarray*}
\begin{eqnarray*}
&=&
\left( 
\frac{1}{n-2} \left( \frac{\kappa }{n-1} + \varepsilon \psi
- \frac{(n^{2} - 3n + 3) \widetilde{\kappa }}{n(n+1)} \right)
+ \frac{ (n-3) \widetilde{\kappa } }{n (n+1)} \right) Q(g,C)\nonumber\\
&=&
\frac{1}{n-2} \left( \frac{\kappa }{n-1} + \varepsilon \psi
- \frac{( 2 n - 3) \widetilde{\kappa }}{n(n+1)} \right)
Q(g,C),
\end{eqnarray*}
i.e. (\ref{DS16Anew01}). Now 
(\ref{02identity01hyper})
and (\ref{DS16Anew01})
lead to (\ref{02identity01hyper17}).
Our proposition is thus proved.
\qed


We finish this section with the following result on quasi-Einstein hypersurfaces 
$M$ in $N_{s}^{n+1}(c)$, $n \geq 4$,
satisfying (\ref{DS4aa}) on ${\mathcal U}_{H} \subset M$, i.e. (\ref{DS4}), with $\rho =0$, 
Precisely, we have  
\begin{thm}
Let $M$ 
be a hypersurface in 
$N_{s}^{n+1}(c)$, $n \geq 4$,
and let (\ref{quasi02}) and (\ref{DS4aa}) be satisfied on ${\mathcal U}_{H} \subset M$.
We have
\newline
(i) The following conditions are satisfied on ${\mathcal U}_{H}$
\begin{eqnarray}
\alpha &=& \frac{\kappa}{n-1} - \frac{ \widetilde{\kappa} }{n(n+1)} ,
\label{qqee02}\\
R \cdot C - C \cdot R &=&  \frac{\kappa}{n-1}\, Q(g,C) - Q(S,C) + A ,
\label{qqee03}
\end{eqnarray}
where the function $\alpha$ is defined by (\ref{quasi02})
and the $(0,6)$-tensor $A$ is defined by  
\begin{eqnarray}
A &=& Q\left(
S -  \frac{1}{n-1} \left( \frac{ (n-2) \kappa }{n-1} +  \frac{ \widetilde{\kappa} }{n(n+1)} \right) g , C \right) .
\label{qqee10}
\end{eqnarray}
(ii) The condition (\ref{Roterformula}) is satisfied on ${\mathcal U}_{H}$ 
if and only if (\ref{qqee08}) holds  on this set. 
\end{thm}
{\bf{Proof.}} (i) If (\ref{DS4aa}) is satisfied on  ${\mathcal U}_{H}$ then (\ref{DZ004Ricciaaa})
holds on this set (see, e.g., Proposition 4.3 (v)). Now (\ref{DZ004Ricciaaa}), 
in view of {\cite[Theorem 2.3] {Ch-DDGP}}, implies 
\begin{eqnarray}
\sum _{ (X_{1}, X_{2}), (X_{3}, X_{4}), (X_{5},X_{6}) } ( R \cdot C - C \cdot R )(X_{1},X_{2},X_{3},X_{4},X_{5},X_{6})
\ =\ 0 .
\label{qqee04}
\end{eqnarray}
Further, (\ref{quasi02}) and  (\ref{qqee04}), in view of {\cite[Proposition 2.1] {DHS105}}, yield (\ref{qqee02}) and
\begin{eqnarray}
(n-2)\, (R \cdot C - C \cdot R) &=& Q(S,R)
- \frac{ \widetilde{\kappa} }{n(n+1)}\, Q(g, R) .
\label{qqee05}
\end{eqnarray}
The condition (\ref{qqee05}), by (\ref{qqee02}), (\ref{eqn2.1}) and (\ref{eqn2.1trtrtr}),  
turns into
\begin{eqnarray}
R \cdot C - C \cdot R &=& \frac{1}{n-2}\, Q(S,C)
- \frac{ \widetilde{\kappa} }{(n-2)n(n+1)}\, Q(g, C) .
\label{qqee06}
\end{eqnarray}
Now from (\ref{qqee06}) we easily get (\ref{qqee03}). 
(ii) From (\ref{qqee03}) it follows immediately that (\ref{Roterformula}) holds on ${\mathcal U}_{H}$ 
if and only if $A = 0$ on this set. Furthermore, 
we note that if $\mathrm{rank}\, (S - \alpha\, g) = 1$ and $\mathrm{rank}\, (S - \beta\, g) = 1$
on ${\mathcal U}_{S} \subset M$, for some functions $\alpha$ and $\beta$, respectively, then  $\alpha = \beta$
at every point of ${\mathcal U}_{S}$ {\cite[Section 3] {DH-2003}}.
Let now $A = 0$ holds on ${\mathcal U}_{H}$. From this, in view of 
{\cite[Proposition 2.4] {DGHHY}} (see also {\cite[Lemma 3.4] {P43}}, 
or, {\cite[Lemma 3.1] {R102}}, or {\cite[Lemma 2.2] {DH-2003}}),
(\ref{quasi02}), (\ref{qqee02})
and the mentioned above remark we easily deduce that
\begin{eqnarray}
\frac{\kappa }{n-1} - \frac{\widetilde{ \kappa }}{n(n+1)} &=& 
\frac{1}{n-1} \left( \frac{ (n-2) \kappa }{n-1} +  \frac{ \widetilde{\kappa} }{n(n+1)} \right)
\label{qqee09}
\end{eqnarray}
on ${\mathcal U}_{H}$. 
From (\ref{qqee09}) we immediately get (\ref{qqee08})(a). 
Now $A = 0$, 
by (\ref{qqee08})(a) and (\ref{qqee09}),
turns into (\ref{qqee08})(b). Conversely, if (\ref{qqee08}) holds on ${\mathcal U}_{H}$ 
then we can check that $A = 0$ on ${\mathcal U}_{H}$.   
\qed

\vspace{2mm}

\noindent
{\bf{Example 4.9.}} 
(i)
As it was mentioned in Proposition 4.3 (i),
for every hypersurface $M$ in $N_{s}^{5}(c)$
(\ref{DS4}) reduces on ${\mathcal U}_{H} \subset M$ to (\ref{DS4aa}).
\newline
(ii)
Let $M$ be a hypersurface in $N_{s}^{n+1}(c)$, $n \geq 4$, satisfying (\ref{DS4aa}) 
on ${\mathcal U}_{H} \subset M$ and
let $Q(S - (\kappa/n)\, g, C) = 0$ (i.e. (\ref{qqee08})(b)) 
at a point $x \in {\mathcal U}_{H}$.  
From the last equation, in view of {\cite[Proposition 2.4] {DGHHY}}
(also see 
{\cite[Proposition 2.1] {DGHZ01}},
or, {\cite[Lemma 3.4] {P43}}, or, {\cite[Lemma 2.2] {DH-2003}}) we have at $x$:
\newline
(a)  
$\mathrm{rank} ( S - (\kappa/n) \, g ) = 1$ and 
\begin{eqnarray*}
\omega (X_{1})\, C(X_{2},X_{3},X_{4},X_{5})
+
\omega (X_{2})\, C(X_{3},X_{1},X_{4},X_{5})
+
\omega (X_{3})\, C(X_{1},X_{2},X_{4},X_{5}) = 1 ,
\end{eqnarray*}
where $\omega$ is some $1$-form at $x$ and
$X_{1}, X_{2}, \ldots , X_{5}$ are vectors tangent to $M$ at $x$, or, 
\newline
(b)  $\mathrm{rank} ( S - (\kappa/n) \, g ) > 1$ and the tensors $C$ and  
$( S - (\kappa/n) \, g ) \wedge ( S - (\kappa/n) \, g )$ are linearly dependent.
From this we get (\ref{eq:h7a}).
An example of a quasi-Einstein hypersurface $M$ in $N_{s}^{n+1}(c)$, $n \geq 4$, 
satisfying 
on ${\mathcal U}_{H} \subset M$ 
the conditions:
(\ref{quasi02}), with $\alpha = \kappa/n$, (\ref{DS4aa}) and (\ref{qqee08})
is given in {\cite[Section 5] {R102}}.
\newline
(iii) 
An example of a non-quasi-Einstein hypersurface $M$ in an Euclidean space ${\mathbb{E}}^{n+1}$, $n \geq 5$,
satisfying (\ref{DS4}) on ${\mathcal U}_{H} = M$, with non-zero functions $\psi $ and $\rho$, is given in
{\cite[Example 5.1 (iii)] {Saw-2015}}. Precisely, on  ${\mathcal U}_{H}$ we have: 
$\rho = ( \kappa\,  \mathrm{tr}(H))/(n-1) \neq 0$, 
$\psi = -  ( \kappa / (n-1) ) $, 
$\widetilde{\kappa } = 0$, 
\begin{eqnarray*}
H^{3} &=& \mathrm{tr}(H)\, H^{2} -  \frac{\kappa}{n-1}\, H +  \frac{\kappa\,  \mathrm{tr}(H) }{n-1}\, g ,
\ \ \
S^{3} \ =\ \frac{2 \kappa }{n-1}\, S^{2} 
+ \left( \frac{ \mathrm{tr}(S^{3}) }{ \kappa } - \frac{2\, \mathrm{tr}(S^{2}) }{n - 1} \right) S  ,
\end{eqnarray*}
\begin{eqnarray}
R \cdot C &=& Q(S,R) + \frac{1}{n-2}\, Q \left( S^{2} - \frac{ \kappa }{n-1}\, S, G \right), 
\label{DZ020}\\
C \cdot R &=& \frac{n-3}{n-2}\, Q(S,R),
\label{DZ021}\\
C \cdot C &=&\frac{n-3}{n-2} \left( Q(S,R) + \frac{1}{n-2}\, Q \left( S^{2} - \frac{\kappa}{n-1}\, S, G \right) \right) . 
\label{DZ022}
\end{eqnarray}
From (\ref{DZ020})-(\ref{DZ022}) it follows immediately that
\begin{eqnarray}
(n-2) ( R \cdot C - C \cdot R) &=&  Q(S,R) + Q \left( S^{2} - \frac{\kappa}{n-1}\, S, G \right) ,
\nonumber\\
C \cdot C &=& C \cdot R + \frac{n-3}{(n-2)^{2}} \, Q \left( S^{2} - \frac{\kappa}{n-1}\, S, G \right) ,
\nonumber\\
C \cdot C   &=& \frac{n-3}{n-2}\, R \cdot C ,
\label{DZ024}
\end{eqnarray}
Using (\ref{DZ024}), $\psi = -  ( \kappa / (n-1) )$, $\varepsilon = 1$ and $\widetilde{\kappa } = 0$ we can check that 
(\ref{DS16Anew01}) holds on $M$. Furthermore, 
using  (\ref{eqn2.1srsr}), (\ref{DZ020}) and (\ref{DZ021}) we obtain 
\begin{eqnarray*}
C \cdot R + \frac{1}{n-2}\, R \cdot C 
&=& 
Q(S,R) - \frac{1}{(n-2)^{2}}\, Q \left( g, g \wedge \left( S^{2} - \frac{\kappa}{n-1}\, S \right) \right)\\
&=&  
Q(S,C) 
- \frac{1}{n-2}\, Q \left( g, \frac{1}{2}\, S \wedge S \right) + \frac{\kappa }{(n-2)(n-1)}\, Q(g, g \wedge S)\\
& &
- \frac{1}{(n-2)^{2}}\, Q \left( g, g \wedge \left( S^{2} - \frac{\kappa}{n-1}\, S \right) \right)\\
&=&
Q(S,C) - \frac{1}{(n-2)^{2}}\, Q \left( g, \frac{n-2}{2}\, S \wedge S - \kappa\, g \wedge S +  g \wedge S^{2} \right) 
\end{eqnarray*}
and in a consequence (\ref{02identity01hyper17}).

\section{The condition $R \cdot C - C \cdot R = L_{1}\, Q(S,C) + L_{2}\, Q(g,C)$}

Let $M$ be a hypersurface in 
a semi-Riemannian space of constant curvature 
$N_{s}^{n+1}(c)$, $n \geq 4$,
satisfying (\ref{cond01})
on ${\mathcal U}_{H} \subset M$,
for some functions 
$L_{1}$ and $L_{2}$ on ${\mathcal U}_{H}$.
We note that in view of {\cite[Corollary 4.1]{2003_DGHV}} 
(\ref{DS4}) holds on ${\mathcal U}_{H}$, and, in a consequence, 
(\ref{ZZ3}) is satisfied on this set.
First we consider quasi-Einstein hypersurfaces satisfying (\ref{cond01}). 
\begin{thm}
If $M$ be a hypersurface in 
$N_{s}^{n+1}(c)$, $n \geq 4$,
satisfying 
on ${\mathcal U}_{H} \subset M$
(\ref{quasi02}) and (\ref{cond01}), for some functions 
$\alpha$, $L_{1}$ and $L_{2}$,
then on this set we have (\ref{DS4aa}), for some function $\psi$, and (\ref{quasi321}),
(\ref{qqee02}), (\ref{qqee03}), (\ref{qqee10}).
Moreover, 
(\ref{Roterformula}) is satisfied on ${\mathcal U}_{H}$ 
if and only if (\ref{qqee08}) holds  on this set. 
\end{thm}
{\bf{Proof.}} From (\ref{cond01}), by an application of Lemma 2.1 (i), we get (\ref{qqee04}).
Now (\ref{quasi02}) and (\ref{qqee04}), in view of {\cite[Proposition 2.1]{DHS105}}, yield (\ref{DZ004Ricciaaa}),
(\ref{qqee02}) and
\begin{eqnarray} 
(n-2)\, ( R \cdot C - C \cdot R ) 
&=& Q(S,R) 
- \frac{\widetilde{\kappa }}{ n (n+1) }\, Q(g,R) .
\label{dhs02}
\end{eqnarray}
From 
{\cite[Proposition 3.2, Theorem 3.1] {1995_DDDVY}}
it follows that on ${\mathcal U}_{H}$ (\ref{DZ004Ricciaaa}) is equivalent to 
(\ref{DS4aa}).
Next, applying to (\ref{dhs02}), the conditions 
(\ref{eqn2.1}),
(\ref{eqn2.1aa})(b)
and
(\ref{eqn2.1trtrtr}),
we get 
(\ref{quasi321}). Now Theorem 4.8 completes the proof.
\qed


In {\cite[Section 4, Example 4.1] {DHS105}} 
an example of a quasi-Einstein non-pseudosym\-met\-ric Ricci-pseudosym\-met\-ric warped product 
$\overline{M} \times _{F} \widetilde{M}$, $\dim \overline{M} = 1$, $\dim \widetilde{N} = n - 1 \geq 4$,
with $\kappa / (n-1) \neq   \widetilde{\kappa } / (n+1)$,
which can be locally realized as a hypersurface $M$ 
in $N_{s}^{n+1}(c)$, $n \geq 5$,
was constructed.
The manifold $(\widetilde{N},\widetilde{g})$ used in that construction is an Einstein manifold. 
Moreover, the condition (\ref{DZ004Ricciaaa}) holds on $M$ {\cite[eq. (4.7)] {DHS105}}.
Thus in view of {\cite[Proposition 3.1 (iii)] {1995_DDDVY}} (\ref{DS4aa}) holds on that hypersurface.
Finally, in view of Theorem 5.1, the equation (\ref{quasi321}) is also satisfied on $M$.


Now we consider non-quasi-Einstein hypersurfaces satisfying (\ref{cond01}). 
Precisely, we consider (\ref{cond01}) at points of ${\mathcal U}_{H}$ 
at which (\ref{dhs_quasi_Einstein}) holds. We have
\begin{thm}
Let $M$ be a hypersurface in 
$N_{s}^{n+1}(c)$, $n \geq 4$, satisfying on ${\mathcal U}_{H} \subset M$
(\ref{cond01}).
If at a point $x \in {\mathcal U}_{H}$ (\ref{dhs_quasi_Einstein})
is satisfied then at this point we have: (\ref{Roterformula}) and 
\begin{eqnarray}
(n-1)\, Q(S,R) 
&=& Q(g,
\left( \varepsilon \psi  + \kappa - \frac{ (n-1) \widetilde{\kappa}}{n (n+1)}   \right) R \nonumber\\
& &
+
\left( \varepsilon \psi - \frac{ 2 (n-1) \widetilde{\kappa} }{n (n+1)} \right) g \wedge S
+ g \wedge S^{2} - \frac{1}{2}\, S \wedge S  )  .
\label{cond02uuu}
\end{eqnarray}
\end{thm} 
{\bf{Proof.}}
The conditions 
(\ref{ZZ3}) and (\ref{cond01}), 
by (\ref{eqn2.1}), (\ref{eqn2.1aa})(a), (\ref{eqn2.1srsr}) and (\ref{DZ005}), 
turn into
\begin{eqnarray}
(n-2)\, (R \cdot C - C \cdot R) &=&
Q(S,R) - \alpha _{3}\, Q(g, g \wedge S) - Q(g , g \wedge S^{2} ) \nonumber\\
& &
+ \left( \frac{(n-1) \widetilde{\kappa} }{n(n+1)} - \frac{\kappa }{n-1} 
- \varepsilon \psi \right)\, Q(g, R),
\label{ZZ3new}
\end{eqnarray}
\begin{eqnarray}
(n-2)\, (R \cdot C - C \cdot R) &=& 
(n-2) L_{1}\, Q(S,R) + L_{1}\, Q(g, \frac{1}{2}\, S \wedge S ) \nonumber\\
& &
+ (n-2) L_{2}\, Q(g,R) - \left( L_{2} + \frac{\kappa L_{1} }{n-1} \right) Q( g, g \wedge S),
\label{cond02}
\end{eqnarray}
respectively.

$I$. The case: $L_{1} = 1 / (n-2)$ at a point $x \in {\mathcal U}_{H}$. 
From (\ref{ZZ3new}) and (\ref{cond02}) we obtain
\begin{eqnarray}
& &
\left(  \frac{(n-1) \widetilde{\kappa} }{n(n+1)} - \frac{\kappa }{n-1} - \varepsilon \psi - (n-2) L_{2} \right) Q(g,R)
- Q(g , g \wedge S^{2})\nonumber\\
& &
- \frac{1}{n-2}\, Q(g, \frac{1}{2}\, S \wedge S )
+ \left( L_{2} + \frac{\kappa }{(n-2)(n-1)}  - \alpha _{3} \right) Q( g, g \wedge S)
\ =\ 0.
\label{1aaa}
\end{eqnarray}

$I(a)$. The subcase: 
\begin{eqnarray*}
L_{2} &=& \frac{1}{n-2} \left( \frac{(n-1) \widetilde{\kappa} }{n(n+1)} - \frac{\kappa }{n-1} - \varepsilon \psi \right) 
\end{eqnarray*}
at a point $x \in {\mathcal U}_{H}$. 
Now (\ref{1aaa}) reduces to
\begin{eqnarray*}
\frac{1}{n-2}\, Q(g, \frac{1}{2}\, S \wedge S ) + Q(g , g \wedge S^{2} )
+ \left( \alpha _{3} -  L_{2} - \frac{\kappa }{(n-2)(n-1)} \right) Q( g, g \wedge S) &=& 0,
\end{eqnarray*}
which yields
\begin{eqnarray*}
\frac{1}{2(n-2)}\, S \wedge S +  g \wedge S^{2} 
+ \left( \alpha _{3} -  L_{2} - \frac{\kappa }{(n-2)(n-1)} \right) g \wedge S + \lambda _{1}\, G &=& 0,
\ \ \  \lambda _{1} \in \mathbb{R},
\end{eqnarray*}
and, by a suitable contraction, leads to
\begin{eqnarray*}
\frac{(n-3)(n-1)}{n-2}\, S^{2} + \left( \alpha _{3} -  L_{2} + \frac{\kappa }{(n-2)(n-1)} \right) S + \lambda _{2}\, g &=& 0, 
\ \ \  \lambda _{2} \in \mathbb{R}.
\end{eqnarray*}
The last two equations yield
$\ (1/2)\, S \wedge S = \beta _{1} \, g \wedge S + \beta _{2}\, G,\ \beta_{1}, \, \beta_{2} \in {\mathbb{R}}$.
From this, in view of {\cite[Lemma 3.1] {2002-G}}, we obtain (\ref{quasi02}), with $\alpha = \beta_{1}$,
which contradicts (\ref{dhs_quasi_Einstein}).
Thus we see that the case $I(a)$ cannot occur at $x$.

$I(b)$. The subcase: 
\begin{eqnarray*}
L_{2} & \neq & \frac{1}{n-2} \left( \frac{(n-1) \widetilde{\kappa} }{n(n+1)} - \frac{\kappa }{n-1} - \varepsilon \psi \right) 
\end{eqnarray*}
at a point $x \in {\mathcal U}_{H}$. Now (\ref{1aaa}) turns into
\begin{eqnarray}
Q(g,R) &=& \alpha_{4}\, Q(g , g \wedge S^{2}) + \frac{\alpha_{4}}{n-2}\, Q\left( g, \frac{1}{2}\, S \wedge S \right)
+ \alpha_{5}\, Q( g, g \wedge S) ,
\label{1aaabbb}\\
\alpha_{4} &=& 
\left(  \frac{(n-1) \widetilde{\kappa} }{n(n+1)} - \frac{\kappa }{n-1} - \varepsilon \psi - (n-2) L_{2} \right)^{-1} ,\nonumber\\
\alpha_{5} &=& \alpha_{4} \left( \alpha _{3}  - L_{2} - \frac{\kappa }{(n-2)(n-1)} \right).\nonumber
\end{eqnarray}
From (\ref{1aaabbb}) we get 
\begin{eqnarray}
R &=& \alpha_{4}\, g \wedge S^{2} + \frac{\alpha_{4}}{2 (n-2) }\,  S \wedge S 
+ \alpha_{5}\, g \wedge S + \lambda_{4} \, G ,\ \ \ \lambda_{4} \in \mathbb{R} .
\label{1aaabbbccc}
\end{eqnarray}
This, by a suitable contraction, yields
\begin{eqnarray}
\frac{(n-3)(n-1)}{n-2}\, S^{2} &=& \alpha_{6}\, S + \lambda_{5} \, g ,\ \ \ \alpha_{6}, \lambda_{5} \in \mathbb{R} .
\label{1aaabbbcccddd}
\end{eqnarray}
From (\ref{1aaabbbccc}) and (\ref{1aaabbbcccddd}) we get
$R = ( \alpha_{4} / (2 (n-2) ))\,  S \wedge S + \alpha_{7}\, g \wedge S + \lambda_{6} \, G$,
$\, a_{7}, \lambda_{6} \in \mathbb{R}$,
which, in wiew of Theorem 2.4, gives (\ref{pseudo}) and (\ref{Roterformula}).
Now $L_{1} = 1 / (n-2)$, (\ref{Roterformula}) and (\ref{cond01}) lead to
\begin{eqnarray*}
Q(S,C) &=& \frac{n-2}{n-1} \left( \frac{ \kappa }{n-1} - L_{2} \right) Q(g,C) ,\\
R \cdot C - C \cdot R &=& \frac{1}{n-1} \left( \frac{ \kappa }{n-1}  + (n-2) L_{2} \right) Q(g,C) .
\end{eqnarray*}
The last condition, in view of {\cite[Theorem 4.1] {DHS-2001}}, yields
\begin{eqnarray}
C \cdot R &=& 0 ,\ \ \ \ R \cdot C \ =\  \frac{1}{n-1} \left( \frac{ \kappa }{n-1}  + (n-2) L_{2} \right) Q(g,C) .
\label{new4abcd}
\end{eqnarray}
As it was stated above, (\ref{pseudo}) holds at $x$. 
Since $x \in {\mathcal U}_{H}$, form Proposition 3.2 (ii) it follows that
$\mathrm{rank}\, H = 2$ and (\ref{pseudo-constant01}) are satisfied at this point.
Further, from Proposition 4.3 (ii) we have (\ref{pseudo-constant02}) and (\ref{pseudo-constant03}). 
Now (\ref{pseudo-constant03}), (\ref{new4abcd}) and $Q(g,R) \neq 0$ imply 
$\widetilde{\kappa}/(n+1) = \kappa /(n-1)$,
which contradicts (\ref{pseudo-constant14}). Thus we see that the case $I(b)$ cannot occur at $x$.

$II$. The case: $L_{1} \neq 1/(n-2)$ at a point $x \in {\mathcal U}_{H}$. 
From (\ref{ZZ3new}) and (\ref{cond02}) it follows that 
the tensor $Q(S,R)$ is a linear combination of the tensors of the form $Q(g,T)$,
where $T$ is a generalized curvature tensor.
Thus, in view of Theorem 4.6, 
then the tensors 
$R \cdot R$, $R \cdot C$, $C \cdot R$, $R \cdot C - C \cdot R$ and $C \cdot C$ satisfy
(\ref{dgpss777}),
(\ref{dgpss201}),
(\ref{dgpss201a}),
(\ref{dgpss201b}) and
(\ref{dgpss777cc}), with the tensors
$B_{1}$, $B_{2}$, $B_{3}$, $B_{4}$ and $B$, 
defined by
(\ref{dgpss26}),
(\ref{dgpss202}),
(\ref{dgpss202a}),
(\ref{dgpss202b})
and
(\ref{dgpss777ccdd}),
respectively.
Now, 
(\ref{dgpss201b}),
(\ref{dgpss202b}) and
(\ref{cond02}) yield
\begin{eqnarray*}
& &
(n-2) \left( - \frac{\varepsilon \psi }{n-1} + \frac{ \widetilde{\kappa}}{n (n+1)} - L_{2} \right) Q(g,R)\nonumber\\
& &
+
\left(
(n-2)\, \left( - \frac{\varepsilon \psi  }{n-1} + \frac{ 2  \widetilde{\kappa} }{n (n+1)} \right)
+ L_{2} + \frac{\kappa L_{1}}{n-1} 
\right) Q( g, g \wedge S)\nonumber\\
& &
- \frac{n-2}{n-1} \, Q(g, g \wedge S^{2} )
- \left( \frac{1}{n-1} +  L_{1} \right) Q\left( g, \frac{1}{2}\, S \wedge S \right) 
\ =\ (n-2) L_{1}\, Q(S,R) 
\end{eqnarray*}
and
\begin{eqnarray}
(n-2) L_{1}\, Q(S,R) 
&=& Q\left( g,
(n-2) \left( - \frac{\varepsilon \psi }{n-1} + \frac{ \widetilde{\kappa}}{n (n+1)} - L_{2} \right) R \right. \nonumber\\
& &
+
\left(
(n-2)\, \left( - \frac{\varepsilon \psi  }{n-1} + \frac{ 2  \widetilde{\kappa} }{n (n+1)} \right)
+ L_{2} + \frac{\kappa L_{1}}{n-1} \right)  g \wedge S\nonumber\\
& &
\left.
- \frac{n-2}{n-1} \,  g \wedge S^{2} 
- \left( \frac{1}{n-1} +  L_{1} \right)  \frac{1}{2}\, S \wedge S  \right)  .
\label{cond02fff}
\end{eqnarray}
But on the other hand,
(\ref{eqn2.1}), (\ref{900ab}) and (\ref{dgpss26}) yield
\begin{eqnarray*}
& &
(n-2) L_{1}\, Q(S,R) \nonumber\\
&=&
Q \left( g, \frac{(n-2)L_{1}}{n-1} 
\left( 
\left( 
\kappa + \varepsilon \psi  - \frac{ (n-1)^{2} \widetilde{\kappa} }{n(n+1)} 
\right) R 
- \frac{1}{2}\, S \wedge S + g \wedge S^{2}
+  
\left( \varepsilon \psi - \frac{ (n - 1) \widetilde{\kappa} }{n(n+1)} \right) g \wedge S \right) \right)\nonumber\\
& &
+ \frac{(n-2)^{2} \widetilde{\kappa} L_{1}}{n (n+1)}\, Q(g,R) 
- \frac{ (n-2) \widetilde{\kappa} L_{1}}{n (n+1)}\, Q(g, g \wedge S) 
\end{eqnarray*}
and
\begin{eqnarray}
(n-2) L_{1}\, Q(S,R) 
&=& Q\left( g,
\frac{(n-2)L_{1}}{n-1}
\left(
\left( \kappa + \varepsilon \psi  - \frac{ (n-1) \widetilde{\kappa} }{n(n+1)} \right) R 
- \frac{1}{2}\, S \wedge S + g \wedge S^{2} \right. \right. \nonumber\\
& &
\left. \left.
+  
\left( \varepsilon \psi - \frac{ 2 (n - 1) \widetilde{\kappa} }{  n(n+1)} \right) g \wedge S \right) \right)  .
\label{cond02ggg}
\end{eqnarray}

Comparing now the right hand sides of 
(\ref{cond02fff}) and (\ref{cond02ggg}) we obtain
\begin{eqnarray}
& & 
\left( - \frac{\varepsilon \psi }{n-1} + \frac{ \widetilde{\kappa}}{n (n+1)} - L_{2}
-
\frac{L_{1}}{n-1}
\left( \kappa + \varepsilon \psi  - \frac{ (n-1) \widetilde{\kappa} }{n(n+1)} \right)
 \right) R\nonumber\\
& &
+
\left(
 - \frac{\varepsilon \psi  }{n-1} + \frac{ 2  \widetilde{\kappa} }{n (n+1)} + \frac{L_{2}}{n-2} 
+ \frac{L_{1}}{n-1}
\left(
\frac{\kappa }{n-2} 
- \varepsilon \psi + \frac{ 2 (n - 1) \widetilde{\kappa} }{ n(n+1)} \right) \right) 
 g \wedge S\nonumber\\
& &
- \frac{1}{n-1} ( 1 + L_{1}) \,  g \wedge S^{2} 
- \frac{1}{(n-2) (n-1) } ( 1 + L_{1}) \frac{1}{2}\,  S \wedge S -   \frac{\lambda_{1}}{n-2} \, G\ =\ 0 . 
\label{5a5a5a}
\end{eqnarray}

The case $II(a)$:  $L_{1} \neq - 1$
at a point $x \in {\mathcal U}_{H}$. From (\ref{5a5a5a}), by a suitable contraction, 
we get
$S^{2} = \beta_{1}\, S + \beta _{2}\, g$, for some $\beta_{1}, \, \beta_{2} \in {\mathbb{R}}$. 
From Proposition 4.3 (ii) it follows that $\rho = 0$, $\beta _{2} = - ( \lambda / n )$ and $\beta _{1} = \alpha _{3}$,
where $\rho $, $\lambda $ and $\alpha _{3}$ are defined by (\ref{DS4}), (\ref{DZ008}) and (\ref{DZ006}), respectively. 
Now (\ref{5a5a5a}) turns into
\begin{eqnarray}
& &
\left( - \frac{\varepsilon \psi }{n-1} + \frac{ \widetilde{\kappa}}{n (n+1)} - L_{2}
-
\frac{L_{1}}{n-1}
\left( \kappa + \varepsilon \psi  - \frac{ (n-1) \widetilde{\kappa} }{n(n+1)} \right)
 \right) R\nonumber\\
&=& 
 \frac{1}{(n-2) (n-1) } ( 1 + L_{1}) \frac{1}{2}\,  S \wedge S
+ \beta _{3} g \wedge S + \beta _{4}\, G ,\ \ \ \beta_{3}, \beta_{4} \in {\mathbb{R}} .  
\label{5a5a5anewbb}
\end{eqnarray}
We note that if 
\begin{eqnarray*}
- \frac{\varepsilon \psi }{n-1} + \frac{ \widetilde{\kappa}}{n (n+1)} - L_{2}
-
\frac{L_{1}}{n-1}
\left( \kappa + \varepsilon \psi  - \frac{ (n-1) \widetilde{\kappa} }{n(n+1)} \right)
&=& 0
\end{eqnarray*}
at a given point then (\ref{5a5a5anewbb}) reduces to 
$(1/2)\, S \wedge S = \beta _{5} \, g \wedge S + \beta _{6}\, G$,
$\beta_{5}, \beta_{6} \in {\mathbb{R}}$.
From this, in view of {\cite[Lemma 3.1] {2002-G}}, we get (\ref{quasi02}), with $\alpha = \beta_{5}$,
which contradicts (\ref{dhs_quasi_Einstein}).
Thus, we see that at this point we must necessary have  
\begin{eqnarray*}
- \frac{\varepsilon \psi }{n-1} + \frac{ \widetilde{\kappa}}{n (n+1)} - L_{2}
-
\frac{L_{1}}{n-1}
\left( \kappa + \varepsilon \psi  - \frac{ (n-1) \widetilde{\kappa} }{n(n+1)} \right)
&\neq & 0 .
\end{eqnarray*}
Now (\ref{5a5a5anewbb}) turns into (\ref{eq:h7a}). Thus, in view of Theorem 2.4 (ii), 
(\ref{Roterformula}) holds at $x$. The conditions (\ref{Roterformula}) and (\ref{cond01}) give 
\begin{eqnarray*}
Q ( (L_{1} + 1)\, S - (L_{2} - \frac{\kappa}{n-1})\, g , C ) &=& 0 .
\end{eqnarray*}
From this,
in view of {\cite[Proposition 2.4] {DGHHY}} (also see {\cite[Proposition 2.1] {DGHZ01}},
or, {\cite[Lemma 3.4] {P43}}, 
or, {\cite[Lemma 2.2] {DH-2003}})
and (\ref{dhs_quasi_Einstein}),
we get
\begin{eqnarray*}
C &=& 
\frac{\lambda }{2} \left( (L_{1} + 1)\, S - \left( L_{2} - \frac{\kappa}{n-1} \right)\, g \right) 
\wedge \left( (L_{1} + 1)\, S - \left( L_{2} - \frac{\kappa}{n-1} \right) g \right), 
\ \ \ \lambda \in {\mathbb{R}}.
\end{eqnarray*}
This, by a suitable contractions, yields 
\begin{eqnarray}
\ \ \
& &
S^{2} \ =\ ( \kappa - (n-2) ( L_{2} - \frac{\kappa }{n-1}) (L_{1} +1 )^{-1})\, S\nonumber\\
\ \ \
& &
- \left( \left( L_{2} - \frac{\kappa }{n-1} \right)^{2} 
+ \left( \kappa (L_{1} + 1) - n \left( L_{2} - \frac{\kappa }{n-1} \right) \right) \left( L_{2} - \frac{\kappa }{n-1} \right) \right)
(L_{1} + 1)^{-2}\, g .
\label{rrr02}
\end{eqnarray}
But on the other hand, in view Theorem 2.4 (i), we have
\begin{eqnarray}
S^{2} &=& 
( \kappa +  ( (n-2) \mu - 1) \phi ^{-1} )\, S +  ( \mu \kappa + (n-1) \eta ) \phi ^{-1} \, g .
\label{rrr03}
\end{eqnarray}
From (\ref{rrr02}) and (\ref{rrr03}) we get
\begin{eqnarray}
 & &
L_{2} - \frac{\kappa }{n-1} \ =\ \frac{1}{n-2} ( 1 - (n-2) \mu ) \phi ^{-1} ) (L_{1} + 1 ) ,
\label{rrr04}\\
& &
\left( 
(n-1) \left( L_{2} - \frac{\kappa }{n-1} \right)
- \kappa (L_{1} + 1)  \right) \left( L_{2} - \frac{\kappa }{n-1} \right)\nonumber\\
& &=\  ( \mu \kappa + (n-1) \eta ) \phi ^{-1}  (L_{1} +1 ) ^{2} .
\label{rrr05}
\end{eqnarray}
Now (\ref{rrr04}) and (\ref{rrr05}) yield
\begin{eqnarray*}
& &
(n-1) ( 1 - (n-2) \mu ) ^{2} - (n-2) \kappa ( 1 - (n-2) \mu ) \phi \ =\ (n-2)^{2} ( \mu \kappa + (n-1) \eta ) \phi ,\\
& &
\kappa \ =\ - (n-2) (n-1) \eta + \frac{n-1}{n-2} ( 1 - (n-2) \mu )^{2} \phi ^{-1} .
\end{eqnarray*}
The last condition, by making use of Proposition 4.3 (iii), leads to 
$\widetilde{\kappa}/(n+1) = \kappa /(n-1)$,
a contradiction with (\ref{pseudo-constant14}). 
Thus we see that the case $II(a)$ cannot occur at $x$.

The case $II(b)$:  $L_{1} = - 1$ at a point $x \in {\mathcal U}_{H}$. 
Thus (\ref{5a5a5a}) reduces to
\begin{eqnarray*}
\left( L_{2} - \frac{\kappa}{n-1} \right) 
\left( R - \frac{1}{n-2}\, g \wedge S \right)  + \lambda _{2}\, G &=& 0 ,\ \ \ \lambda_{2} \in {\mathbb{R}} ,
\end{eqnarray*}
which yields
\begin{eqnarray*}
\left( L_{2} -  \frac{ \kappa}{n-1}  \right) \left( R - \frac{1}{n-2}\, g \wedge S + \frac{\kappa }{(n-2)(n-1)}\, G \right)   
&=&  \lambda _{3}\, G ,\ \ \ \lambda_{3} \in {\mathbb{R}} .
\end{eqnarray*}
This, by (\ref{eqn2.1}), turns into
$( L_{2} - \kappa /(n-1) ) \, C   = \lambda _{3}\, G$, 
which implies $( L_{2} - \kappa /(n-1) ) \, C  = 0$ and, in a consequence, $L_{2} = \kappa /(n-1)$. 
Now (\ref{cond02fff}) reduces to (\ref{cond02uuu}). 
Our theorem is thus proved.
\qed


As a consequence of 
{\cite[Corollary 4.1] {2003_DGHV}} (see also Section 4) and theorems 4.6 and 5.2 we get

\begin{thm}
Let $M$ be a hypersurface in 
$N_{s}^{n+1}(c)$, $n \geq 4$, 
and let (\ref{dhs_quasi_Einstein}) and (\ref{cond01}) 
be satisfied at every point of the set ${\mathcal U}_{H} \subset M$. 
Then (\ref{Roterformula}) and (\ref{cond02uuu}) hold on ${\mathcal U}_{H}$.
Moreover, the tensors 
$R \cdot R$, $R \cdot C$, $C \cdot R$, $R \cdot C - C \cdot R$ and $C \cdot C$ satisfy
(\ref{dgpss777}), (\ref{dgpss201}), (\ref{dgpss201a}), (\ref{dgpss201b}) and (\ref{dgpss777cc}), 
with the tensors $B_{1}$, $B_{2}$, $B_{3}$, $B_{4}$ and $B$, 
defined by (\ref{dgpss26}), (\ref{dgpss202}), (\ref{dgpss202a}), (\ref{dgpss202b}) and (\ref{dgpss777ccdd}), respectively. 
\end{thm}


\noindent
{\bf{Acknowledgments}}

The first two named authors of this paper are supported 
by the Wroc\l aw University of Environmental and Life Sciences, Poland.

\vspace{5mm}

\noindent
\footnotesize{Ryszard Deszcz and Ma\l gorzata G\l ogowska\\
Wroc\l aw University of Environmental and Life Sciences\\ 
Department of Mathematics\\
Grunwaldzka 53, 50-357 Wroc\l aw, Poland}\\
{\sf E-mail: ryszard.deszcz@upwr.edu.pl},
{\sf E-mail: malgorzata.glogowska@upwr.edu.pl}

\vspace{1mm}

\noindent
\footnotesize{Georges Zafindratafa\\
Universit\'{e} Polytechnique des Hauts de France\\
Laboratoire de Math\'{e}matiques et Applications de Valenciennes\\
Le Mont Houy, LAMAV-ISTV2, 59313 Valenciennes Cedex 9, France}\\
{\sf E-mail: georges.zafindratafa@uphf.fr}


\begin{thebibliography}{99}

\footnotesize

\bibitem{AD-2002} 
B.E. Abdalla and F. Dillen, 
\emph{A Ricci-semi-symmetric hypersurface of the Euclidean
space which is not semi-symmetric}, 
Proc. Amer. Math. Soc. ${\mathbf{130}}$ (2002), 1805--1808.


\bibitem{ARS-1995} L. Al\' \i as, A. Romero and M. S\'{a}nchez,
\emph{Compact spacelike hypersurfaces of constant mean curvature in generalized
Robertson-Walker spacetimes}, 
in: Geometry and Topology of Submanifolds, 
$\mathbf{VII}$, World Sci., River Edge, NJ, 1995, 67--70. 


\bibitem{ARS-1997} L. Al\' \i as, A. Romero and M. S\'{a}nchez,
\emph{Spacelike hypersurfaces of constant mean curvature and Calabi-Einstein type
problems}, T\^{o}hoku Math. J. ${\mathbf{49}}$ (1997), 337--345.


\bibitem{ADEMO-2002} 
K. Arslan, R. Deszcz, R. Ezenta\c s, C. Murathan, and C. \"{O}zg\"{u}r,
\emph{On pseudosymmetry type hypersurfaces of semi-Euclidean spaces. I}, 
Acta Math. Sci. ${\mathbf{22 B}}$ (2002), 346--358.


\bibitem{P119} 
K. Arslan, R. Deszcz, R. Ezenta\c s, M. Hotlo\'{s} and C. Murathan,
\emph{On generalized Robertson-Walker spacetimes satisfying some curvature condition},
Turkish J. Math. ${\mathbf{38}}$ (2014), 353--373. 


\bibitem{ADY-1997} 
K. Arslan, R. Deszcz and \c S. Yaprak, 
\emph{On Weyl pseudosymmetric hypersurfaces}, 
Colloq. Math. ${\mathbf{72}}$ (1997), 353--361.

\bibitem{Besse-1987} 
A.L. Besse, 
\emph{Einstein Manifolds},
Ergeb. Math. Grenzgeb., 3. Folge, Bd. 10. Springer-Verlag, Berlin, Heidelberg, New York, 1987.


\bibitem{TEC_GRJ_1998} 
T.E. Cecil and G.R. Jensen, 
\emph{Dupin hypersurfaces with three principal curvatures}, 
Invent. Math. ${\mathbf{132}}$ (1998), 121--178.


\bibitem{TEC_PJR_2015}
T.E. Cecil and P.J. Ryan, 
\emph{Geometry of Hypersurfaces}, Springer Monographs in Mathematics, 
Springer New York Heidelberg Dodrecht London, 2015.


\bibitem{Chen-2011}
B.-Y. Chen, 
\emph{Pseudo-Riemannian Geometry}, $\delta$\emph{-Invariants and Applications},
World Scientific, 2011.


\bibitem{Chen-2017}
B.-Y. Chen, 
\emph{Differential Geometry of Warped Product Manifolds and Submanifolds},
World Scientific, 2017.


\bibitem{Chen-2017-KJM}
B.-Y. Chen, 
\emph{Classification of torqued vector fields and its applications to Ricci solitons},
Kragujevac J. Math. ${\mathbf{41}}$ (2017), 239--250.


\bibitem{Chen-Yildirim-2015}
B.-Y. Chen and H. Y{\i}ld{\i}r{\i}m, 
\emph{Classification of ideal submanifolds of real space forms with type number $\leq 2$},
J. Geom. Phys. ${\mathbf{92}}$ (2015), 167--180.


\bibitem{Ch-DDGP} 
J. Chojnacka-Dulas, R. Deszcz, M. G\l ogowska and M. Prvanovi\'{c}, 
\emph{On warped products manifolds satisfying some curvature conditions},
J. Geom. Phys. ${\mathbf{74}}$ (2013), 328--341.


\bibitem{DecuP-TSVer}
S. Decu, M. Petrovi\'{c}-Torga\v{s}ev, A. \v{S}ebekovi\'{c} and L. Verstraelen,
\emph{On the Roter type of Wintgen ideal submanifolds},
Rev. Roumaine Math. Pures Appl. ${\mathbf{57}}$ (2012), 75--90. 


\bibitem{1995_DDDVY} 
F. Defever, R. Deszcz, P. Dhooghe, L. Ver\-stra\-elen and \c S. Yaprak, 
\emph{On Ricci-pseudo\-sym\-met\-ric hyper\-sur\-faces in spaces of constant curvature}, 
Results in Math. ${\mathbf{27}}$ (1995), 227--236.


\bibitem{2000_DDKV} F. Defever, R. Deszcz and D. Kowalczyk,
\emph{Semisymmetry and Ricci-semisymmetry for hypersurfaces of
semi-Riemannian space form}, 
Arab J. Math. Sc. ${\mathbf{6}}$ (2000), 1--16.


\bibitem{49} 
F. Defever, R. Deszcz and M. Prvanovi\'{c}, 
\emph{On warped product manifolds satisfying some curvature condition 
of pseudosymmetry type},
Bull. Greek Math. Soc. ${\mathbf{36}}$ (1994), 43--62.


\bibitem{1997_DDSVY} 
F. Defever, R. Deszcz, Z. \c Sent\"{u}rk, L. Verstraelen and \c S. Yaprak, 
\emph{On a problem of P.J. Ryan},
Kyungpook Math. J. 
${\mathbf{37}}$ (1997), 371--376.


\bibitem{1999_DDSVY} 
F. Defever, R. Deszcz, Z. \c Sent\"{u}rk, L. Verstraelen and \c S. Yaprak, 
\emph{P.J. Ryan's problem in semi-Riemannian space forms},
Glasgow Math. J. ${\mathbf{41}}$ (1999), 271--281.


\bibitem{1989_DDV} J. Deprez, R. Deszcz and L. Verstraelen, 
\emph{Examples of pseudosymmetric conformally flat warped products}, 
Chinese J. Math. ${\mathbf{17}}$ (1989), 51--65.


\bibitem{D-1992} R. Deszcz, 
\emph{On pseudosymmetric spaces}, 
Bull. Soc. Math. Belg.${\mathbf{44}}$ (1992), S\'{e}r. A, Fasc. 1, 1--34.

\bibitem{D-1997} R. Deszcz, 
\emph{On pseudosymmetric hypersurfaces in spaces of constant curvature},
Tensor (N.S.) ${\mathbf{58}}$ (1997), 253--269. 

\bibitem{P106} 
R. Deszcz, 
\emph{On some Akivis-Goldberg type metrics},
Publ. Inst. Math. (Beograd) (N.S.)
$\mathbf{74}$ ($\mathbf{88}$) (2003), 71--83.


\bibitem{2002-DG-1} R. Deszcz and M. G\l ogowska,
\emph{Examples of nonsemisymmetric Ricci-semisymmetric hypersurfaces},
Colloq. Math. ${\mathbf{94}}$ (2002), 87--101. 


\bibitem{DG90} 
R. Deszcz and M. G\l ogowska, 
\emph{Some nonsemisymmetric Ricci-semisymmetric
warped product hypersurfaces},
Publ. Inst. Math. (Beograd) (N.S.) 
$\mathbf{72}$ ($\mathbf{86}$) (2002), 81--94.


\bibitem{DGHHY} 
R. Deszcz, M. G\l ogowska, H. Hashiguchi, M. Hotlo\'{s} and M. Yawata,
\emph{On semi-Riemannian manifolds satisfying some conformally invariant curvature condition}, 
Colloq. Math. ${\mathbf{131}}$ (2013), 149--170.


\bibitem{DGHJZ01}
R. Deszcz, M. G\l ogowska, M. Hotlo\'{s}, J. Je\l owicki and G. Zafindratafa,
\emph{Curvature properties of some warped product manifolds},
poster,
conference Differential Geometry, Banach Conference Center at Bedlewo, Poland,
June 19 - June 24, 2017.


\bibitem{DGHSaw}  
R. Deszcz, M. G\l ogowska, M. Hotlo\'{s}, and K. Sawicz,
\emph{A Survey on Generalized Einstein Metric Conditions},
Advances in Lorentzian Geometry:
Proceedings of the Lorentzian Geometry Conference in Berlin,
AMS/IP Studies in Advanced Mathematics $\mathbf{49}$, S.-T. Yau (series ed.),
M. Plaue, A.D. Rendall and M. Scherfner (eds.), 2011, 27--45.


\bibitem{DGHS} 
R. Deszcz, M. G\l ogowska, M. Hotlo\'{s} and Z. \c Sent\"{u}rk,
\emph{On certain quasi-Einstein semisymmetric hypersurfaces},
Ann. Univ. Sci. Budapest. E\"{o}tv\"{o}s Sect. Math. 
${\mathbf{41}}$ (1998), 151--164.


\bibitem{2003_DGHV} 
R. Deszcz, M. G\l ogowska, M. Hotlo\'{s}, and L. Verstraelen,
\emph{On some generalized Einstein metric conditions 
on hypersurfaces in semi-Riemannian space forms},
Colloq. Math. ${\mathbf{96}}$ (2003), 149--165.


\bibitem{DGHZ01}
R. Deszcz, M. G\l ogowska, M. Hotlo\'{s} and G. Zafindratafa,
\emph{On some curvature conditions of pseudosymmetry type}, 
Period. Math. Hung. ${\mathbf{70}}$ (2015), 153--170.


\bibitem{2016_DGHZhyper}
R. Deszcz, M. G\l ogowska, M. Hotlo\'{s} and G. Zafindratafa,
\emph{Hypersurfaces in spaces forms satisfying 
some curvature conditions},
J. Geom. Phys. ${\mathbf{99}}$ (2016), 218--231.


\bibitem{DGJZ-2016} 
R. Deszcz, M. G\l ogowska, J. Je\l owicki and G. Zafindratafa,
\emph{Curvature properties of some class of warped product manifolds}, 
Int. J. Geom. Meth. Modern Phys. ${\mathbf{13}}$ (2016), 1550135 (36 pages).


\bibitem{DGP-TV01}
R. Deszcz, M. G\l ogowska, M. Petrovi\'{c}-Torga\v{s}ev and L. Verstraelen,
\emph{On the Roter type of Chen ideal submanifolds},
Results. Math. ${\mathbf{59}}$ (2011), 401--413.


\bibitem{DGP-TV02}
R. Deszcz, M. G\l ogowska, M. Petrovi\'{c}-Torga\v{s}ev and L. Verstraelen,
\emph{Curvature properties of some class of minimal hypersurfaces in Euclidean spaces},
Filomat ${\mathbf{29}}$ (2015), 479--492.


\bibitem{2011-DGPSS} 
R. Deszcz, M. G\l ogowska, M. Plaue, K. Sawicz and M. Scherfner, 
\emph{On hypersurfaces in space forms satisfying particular curvature conditions of Tachibana type}, 
Kragujevac J. Math. ${\mathbf{35}}$ (2011), 223--247.


\bibitem{DHV2008}
R. Deszcz, S. Haesen and L. Verstraelen,  
\emph{On natural symmetries}, Topics in Differential Geometry, Ch. 6, 
Editors A. Mihai, I. Mihai and R. Miron, 
Editura Academiei Rom$\hat{\mbox{a}}$ne, 2008, 249--308.


\bibitem{P43} R. Deszcz and M. Hotlo\'{s}, 
\emph{On certain subclass of pseudosymmetric manifolds}, 
Publ. Math. Debrecen $\mathbf{53}$ (1998), 29--48.


\bibitem{R102}
R. Deszcz and M. Hotlo\'{s},
\emph{On hypersurfaces with type number two in space forms},
Ann. Univ. Sci. Budapest. E\"{o}tv\"{o}s Sect. Math. 
$\mathbf{46}$ (2003), 19--34.


\bibitem{DH-2003}
R. Deszcz and M. Hotlo\'{s}, 
\emph{On some pseudosymmetry type curvature condition},
Tsukuba J. Math. $\mathbf{27}$  (2003) 13--30.


\bibitem{DH-2018}
R. Deszcz and M. Hotlo\'{s}, 
\emph{On geodesic mappings in particular class of Roter spaces},
arXiv: 1812.00670 [math.DG], submitted 3 December 2018.  


\bibitem{DeHoJJKunSh}
R. Deszcz, M. Hotlo\'{s}, J. Je\l owicki, H. Kundu, and A.A. Shaikh,
\emph{Curvature properties of G\"{o}del metric}.
Int. J. Geom. Meth. Modern Phys. $\mathbf{11}$ (2014), 1450025 (20 pages).


\bibitem{DHS-2001}
R. Deszcz, M. Hotlo\'{s} and Z. \c Sent\"{u}rk,  
\emph{On some family of generalized Einstein
metric curvature conditions}, 
Demonstratio Math. $\mathbf{34}$ (2001), 943--954.


\bibitem{DHS105}
R. Deszcz, M. Hotlo\'{s} and Z. \c Sent\"{u}rk,
\emph{On curvature properties of certain quasi-Einstein hypersurfaces},
Int. J. Math. $\mathbf{23}$ (2012), 1250073 (17 pages).


\bibitem{DeKow}
R. Deszcz and D. Kowalczyk,
\emph{On some class of pseudosymmetric warped products},
Colloq. Math.
$\mathbf{97}$ (2003), 7--22.


\bibitem{DP-TVZ} 
R. Deszcz, M. Petrovi\'{c}-Torga\v{s}ev, L. Verstraelen and G. Zafindratafa,
\emph{On Chen ideal submanifolds satisfying some conditions of pseudo-symmetry type},
Bull. Malaysian Math. Sci. Soc. $\mathbf{39}$ (2016), 103--131.


\bibitem{DePlaScher}  
R. Deszcz, M. Plaue and M. Scherfner, 
\emph{On Roter type warped products with 1-dimensional fibres}, 
J. Geom. Phys. $\mathbf{69}$ (2013), 1--11.


\bibitem{DeScher}
R. Deszcz and M. Scherfner, 
\emph{On a particular class of warped products with fibres locally
isometric to generalized Cartan hypersurfaces},
Colloq. Math. $\mathbf{109}$ (2007), 13--29. 


\bibitem{2005-DSaw} 
R. Deszcz and K. Sawicz, 
\emph{On some class of hypersurfaces in Euclidean spaces},
Annales Univ. Sci. Budapest. E\"{o}tv\"{o}s Sect. Math.
${\mathbf{48}}$ (2005), 87--98.


\bibitem{DV-1991} 
R. Deszcz and L. Verstraelen, 
\emph{Hypersurfaces of semi-Riemannian conformally flat manifolds}, 
in: Geometry and Topology of Submanifolds, ${\mathbf{III}}$,
World Sci., River Edge, NJ, 1991, 131--147. 


\bibitem{DVV1991} 
R. Deszcz, L. Verstraelen and L. Vrancken, 
\emph{The symmetry of warped product spacetimes}, 
Gen. Relativ. Gravit. $\mathbf{23}$ (1991), 671--681.


\bibitem{DVY} 
R. Deszcz, L. Verstraelen and \c S. Yaprak, 
\emph{Pseudosymmetric hypersurfaces in 4-dimensional spaces of constant 
curvature}, 
Bull. Inst. Math. Acad. Sinica ${\mathbf{22}}$ (1994), 167--179.


\bibitem{DY-1994-cm} 
R. Deszcz and \c S. Yaprak, 
\emph{Curvature properties of Cartan hypersurfaces}, 
Colloq. Math. $\mathbf{67}$ (1994), 91--98.


\bibitem{DY1994} 
R. Deszcz and \c S. Yaprak, 
\emph{Curvature properties of certain pseudosymmetric manifolds}, 
Publ. Math. Debrecen $\mathbf{45}$ (1994), 333--345.


\bibitem{DVZ-2018}
M.M. Diniz, J.A.M. Vilhena and J.F.Z. Zapata,
\emph{Minimal hypersurfaces in} $\mathbb{S}^{5}$ \emph{with vanishing Gauss-Kronecker curvature},
Izv. Math. {\bf{82}} (2018), 477--493.


\bibitem{2002-G}
M. G\l ogowska,
\emph{Semi-Riemannian manifolds whose Weyl tensor is a Kulkarni-Nomizu square},
Publ. Inst. Math. (Beograd) (N.S.) 
$\mathbf{72}$ ($\mathbf{86}$) (2002), 95--106.


\bibitem{2005-G}
M. G\l ogowska, 
\emph{Curvature conditions on hypersurfaces with two distinct principal curvatures},
in: Banach Center Publ. ${\mathbf{69}}$, Inst. Math. Polish Acad. Sci., 2005, 133--143.


\bibitem{G5} 
M. G\l ogowska,
\emph{On Roter type manifolds},
in: Pure and Applied Differential Geometry - PADGE 2007, Shaker Verlag, Aachen, 2007, 114--122.


\bibitem{G6}
M. G\l ogowska,
\emph{On quasi-Einstein Cartan type hypersurfaces},
J. Geom. Phys. $\mathbf{58}$ (2008), 599--614.


\bibitem{GrifPod}
J.B. Griffiths and J. Podolsk\'{y},
\emph{Exact Space-Times in Einstein's General Relativity},
Cambridge Univ. Press, 2009.


\bibitem{HV_2007} 
S. Haesen and L. Verstraelen,
{\em Properties of a scalar curvature invariant depending on two planes},
Manuscripta Math. $\mathbf{122}$ (2007), 59--72.


\bibitem{HaVerSigma}
S. Haesen and L. Verstraelen,  
\emph{Natural intrinsic geometrical symmetries},
SIGMA ${\mathbf{5}}$ (2009), 086 (14 pages).


\bibitem{Kow01}  
D. Kowalczyk, 
\emph{On some class of semisymmetric manifolds}, 
Soochow J. Math. $\mathbf{27}$ (2001), 445--461.


\bibitem{Kow02} 
D. Kowalczyk,
\emph{On the Reissner-Nordstr\"{o}m-de Sitter type spacetimes},
Tsukuba J. Math. ${\mathbf{30}}$ (2006), 363--381.  


\bibitem{KowSek_1997}
O. Kowalski and M. Sekizawa,
{\em Pseudo-symmetric spaces of constant type
in dimension three -- elliptic spaces},
Rend. Mat. Appl., VII Ser.,
$\mathbf{17}$ (1997), 477--512.


\bibitem{Lumiste} 
\"{U}. Lumiste, 
\emph{Semiparallel Submanifolds in Space Forms}, Springer Science $+$ Business Media,
New York, LLC 2009.


\bibitem{MSV_2015} J. Mike\v{s}, E. Stepanova, A. Van\v{z}urova and et al.,
{\sl{Differential geometry of special mappings}},
Palack\'{y} Univ. Olomouc, Fac. Sci., Olomouc, 2015.  


\bibitem{Mir-1992}
V.A. Mirzoyan, 
\emph{Structure theorems for Riemannian Ric-semisymmetric spaces}, 
Izv. Vyssh. Uchebn. Zaved. Mat.
$\mathbf{6}$ (1992), 80--89 (in Russian); English transl.: 
Russ. Math. $\mathbf{36}$ ($\mathbf{6}$) (1992), 75--83.


\bibitem{Mir-1999}
V.A. Mirzoyan, 
\emph{Submanifolds with parallel and semi-parallel structures},
J. Contemp. Math. Anal. Armenian Acad. Sci. $\mathbf{34}$ (1999), 69--73. 


\bibitem{Mir-Mach-2012} 
V.A. Mirzoyan and G.S. Machkalyan, 
\emph{Normally flat Ric-semisymmetric submanifolds in Euclidean spaces},
Izv. Vyssh. Uchebn. Zaved. Mat. $\mathbf{9}$ (2012), 19--31 (in Russian); English transl.:
Russ. Math. $\mathbf{56}$ ($\mathbf{9}$) (2012), 14--24.


\bibitem{Ryan-1972}
P.J. Ryan, 
\emph{A class of complex hypersurfaces}, 
Colloq. Math. $\mathbf{26}$ (1972), 175--182.


\bibitem{Saw-2004} 
K. Sawicz, 
\emph{Hypersurfaces in spaces of constant curvature
satisfying some Ricci-type conditions}, 
Colloq. Math. $\mathbf{101}$ (2004), 183--201.

\bibitem{Saw-2005} 
K. Sawicz,
\emph{Curvature conditions on hypersurfaces with three distinct principal curvatures},
in: Banach Center Publications ${\mathbf{69}}$, Inst. Math. Polish Acad. Sci., 2005, 145--156.

\bibitem{Saw-2006}
K. Sawicz, 
\emph{On curvature characterization of
some hypersurfaces in spaces of constant curvature}, 
Publ. Inst. Math. (Beograd) (N.S.)
${\mathbf{79}}$ (${\mathbf{93}}$) (2006), 95--107.

\bibitem{Saw-2007}
K. Sawicz,
\emph{Curvature identities on hypersurfaces in semi-Riemannian space forms},
in: Pure and Applied Differential Geometry, PADGE 2007, Shaker Verlag,
Aachen 2007, 114--122.

\bibitem{Saw-2015}
K. Sawicz,
\emph{Curvature properties of some class of hypersurfaces in Euclidean spaces}, 
Publ. Inst. Math. (Beograd) (N.S.) ${\mathbf{98}}$ (${\mathbf{112}}$) (2015), 165--177.


\bibitem{SDHJK} 
A.A. Shaikh, R. Deszcz, M. Hotlo\'{s}, J. Je\l owicki, and H. Kundu,
\emph{On pseudosymmetric manifolds}, 
Publ. Math. Debrecen $\mathbf{86}$ (2015), 433--456. 


\bibitem{SH}
Z. Stuchlik and S. Hledik,
\emph{Properties of the Reissner-Nordstr\"{o}m spacetimes with a nonzero cosmological constant},
Acta Phys. Slovaca $\mathbf{52}$ (2002), 363--407. 


\bibitem{Szabo}
Z.I. Szab\'{o}, 
\emph{Structure theorems on Riemannian spaces satisfying} $R(X,Y) \cdot R = 0$. 
\emph{I. The local version}, J. Differential Geom. $\mathbf{17}$ (1982), 531--582.


\bibitem{LV1}
L. Verstraelen, 
\emph{Comments on the pseudo-symmetry in the sense of Ryszard Deszcz}, in: 
Geometry and Topology of Submanifolds, $\mathbf{VI}$, World Sci., Singapore, 1994,
119--209.


\bibitem{LV2} 
L. Verstraelen, 
\emph{Natural extrinsic geometrical symmetries - an introduction -},
in: Recent Advances in the Geometry of Submanifolds 
Dedicated to the Memory of Franki Dillen (1963-2013), 
Contemporary Mathematics, $\mathbf{674}$ (2016), 5--16.


\bibitem{LV3-Foreword}
L. Verstraelen, 
\emph{Foreword}, in:  B.-Y. Chen,
\emph{Differential Geometry of Warped Product Manifolds and Submanifolds},
World Scientific, 2017, vii--xxi.
\end{thebibliography}
\end{document}